\documentclass[a4paper,12pt]{article}     

\usepackage{graphicx}		  
\usepackage{amsmath,amssymb}	

\newtheorem{lemma}{Lemma}[section]
\newtheorem{theorem}[lemma]{Theorem}
\newtheorem{proposition}[lemma]{Proposition}
\newtheorem{corollary}[lemma]{Corollary}
\newtheorem{assumption}{Assumption}

\def\authorfont{\footnotesize}

\def\ccode#1{\par
\vspace*{8pt}
{\authorfont{\leftskip18pt\rightskip\leftskip
\noindent #1\par}}\par}

\newenvironment{Proof}{
\hspace*{-9mm}
{ \it Proof.}}
{\hfill {$\square$}\vspace{1.5em}}

\begin{document}

\begin{center}{
{\Large 
 Properties of minimal charts and
 their applications X: 
 charts of type $(5,2)$}
\vspace{10pt}
\\ 
Teruo NAGASE and Akiko SHIMA\footnote{The second author is supported by JSPS KAKENHI Grant Number 21K03255.}
}
\end{center}

\begin{abstract}
Charts are oriented labeled graphs in a disk.
Any simple surface braid (2-dimensional braid) can be described by using a chart.
Also, a chart represents an oriented closed surface
embedded in 4-space.
In this paper, we investigate embedded surfaces in 4-space
by using charts.
Let $\Gamma$ be a chart,
and we denote by $\Gamma_m$
the union of all the edges of label $m$.
A chart $\Gamma$ is of type $(5,2)$
if there exists a label $m$
such that 
$w(\Gamma)=7$,
$w(\Gamma_m\cap\Gamma_{m+1})=5$,
$w(\Gamma_{m+1}\cap\Gamma_{m+2})=2$
where 
$w(G)$ is the number of white vertices in $G$.
In this paper, we investigate
a minimal chart of 
type $(5,2)$.
\end{abstract}

%
%
%
%

\ccode{2020 Mathematics Subject Classification. Primary 57K45,05C10; Secondary 57M15.}
\ccode{ {\it Key Words and Phrases}. surface link, chart, C-move, white vertex. }


\setcounter{section}{0}
\section{Introduction}


Charts are oriented labeled graphs in a disk (see  \cite{KnottedSurfaces},\cite{BraidBook}, and see Section~\ref{s:Prel}  for the precise definition of charts).
Let $D_1^2, D_2^2$ be 2-dimensional disks.
Any simple surface braid (2-dimensional braid) can be described 
by using a chart,
here a simple surface braid is a properly embedded surface
$S$ in the 4-dimensional disk $D_1^2\times D_2^2$ such that
a natural map $\pi:S\subset D_1^2\times D_2^2\to D_2^2$ is 
a simple branched covering map of $D_2^2$ and
the boundary $\partial S$ is a trivial closed braid in
the solid torus $D_1^2\times \partial D_2^2$
(see \cite{BraidThree}, \cite[Chapter 14 and Chapter 18]{BraidBook}).
Also, from a chart, 
we can construct a simple closed surface braid in 4-space ${\Bbb R}^4$. This surface is an oriented closed surface 
embedded in ${\Bbb R}^4$.
On the other hand, any oriented embedded closed surface 
 in ${\Bbb R}^4$ is ambient isotopic to a simple
closed surface braid
 (see \cite{BraidThree},\cite[Chapter 23]{BraidBook}). 
A C-move 
is a local modification between two charts
in a disk (see Section~\ref{s:Prel} for C-moves).
A C-move between two charts induces 
an ambient isotopy between oriented closed surfaces 
corresponding to the two charts.
In this paper, we investigate oriented closed surfaces in 4-space
by using charts.

We will work in the PL category or smooth category. All submanifolds are assumed to be locally flat.
In \cite{ONS},
we showed that there is no minimal chart with exactly five vertices
 (see Section~\ref{s:Prel} for the precise definition of minimal charts). 
Hasegawa proved that there exists a minimal chart with exactly
six white vertices \cite{H1}. 
This chart represents a 2-twist spun trefoil.
In \cite{INS} and \cite{NST},
we investigated minimal charts with exactly four white vertices.
In this paper, 
we investigate properties of minimal charts 
which support a conjecture that
there is no minimal chart with exactly seven white vertices
(see \cite{ChartApp1},\cite{ChartAppII},\cite{ChartAppIII},\cite{ChartAppIV},\cite{ChartAppV},\cite{ChartAppVI},\cite{ChartAppVII},\cite{ChartAppVIII},\cite{ChartAppIX},\cite{ChartAppXI}).

Let $\Gamma$ be a chart.
For each label $m$, we denote by $\Gamma_m$
the union of all the edges of label $m$.

Now we define a type of a chart:
Let $\Gamma$ be a chart with at least one white vertex, 
and $n_1,n_2,\dots,n_k$ integers.
The chart $\Gamma$ is of {\it type $(n_1,n_2,\dots,n_k)$} if there exists a label $m$ of $\Gamma$ satisfying the following three conditions:
\begin{enumerate}
\item[(i)] For each $i=1,2,\dots, k$, 
the chart $\Gamma$ contains exactly $n_{i}$ white vertices in $\Gamma_{m+i-1}\cap \Gamma_{m+i}$.
\item[(ii)] If $i<0$ or $i>k$, then $\Gamma_{m+i}$ does not contain any white vertices.
\item[(iii)] Both of the two subgraphs $\Gamma_m$ and $\Gamma_{m+k}$ contain at least one white vertex.
\end{enumerate}
If we want to emphasize the label $m$,
then we say that $\Gamma$ is of {\it type $(m;n_1,n_2,\dots,n_k)$}. 
Note that $n_1\ge1$ and $n_k\ge1$ by Condition~(iii).

We proved in \cite[Theorem 1.1]{ChartAppII} that
if there exists a minimal $n$-chart $\Gamma$ with exactly seven white vertices,
then $\Gamma$ is a chart of 
type $(7),(5,2),(4,3),(3,2,2)$ or $(2,3,2)$ 
(if necessary we change the label
$i$ by $n-i$ for all label $i$).
In \cite{ChartAppV},
we showed that
there is no minimal chart of type $(3,2,2)$.
In \cite{ChartAppVI} and \cite{ChartAppVII},
there is no minimal chart of type $(2,3,2)$.
In \cite{ChartAppVIII},
there is no minimal chart of type $(7)$.
In \cite{ChartAppIX},
there is no minimal chart of type $(4,3)$.
In this paper, we investigate a minimal chart of type $(5,2)$.

An edge in a chart is called 
a {\it terminal edge}
if it has
a white vertex and a black vertex.

\begin{figure}[htb]
\centerline{\includegraphics{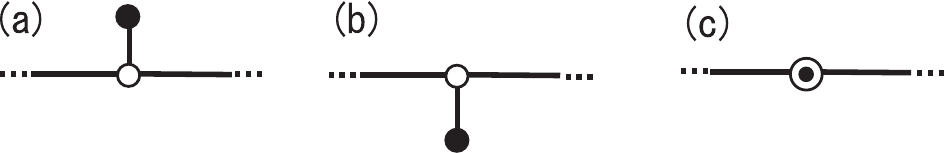}}
\caption{\label{Fig01}
(a),(b) White vertices in terminal edges.
(c) BW-vertex.}
\end{figure}

 In our argument  we often construct a chart $\Gamma$. 
On the construction of a chart $\Gamma$, for a white vertex $w\in\Gamma_m$ for some label $m$,  
among the three edges of $\Gamma_m$ 
containing $w$, 
if one of the three edges is a terminal edge 
(see Fig.~\ref{Fig01}(a) and (b)), 
then we remove the terminal edge and
put a black dot at the center of the white vertex  as shown in Fig.~\ref{Fig01}(c).
Namely
Fig.~\ref{Fig01}(c) means 
Fig.~\ref{Fig01}(a) or 
Fig.~\ref{Fig01}(b).
We call the vertex in Fig.~\ref{Fig01}(c) 
a {\it BW-vertex} with respect to $\Gamma_m$.

\begin{figure}[htb]
\centerline{\includegraphics{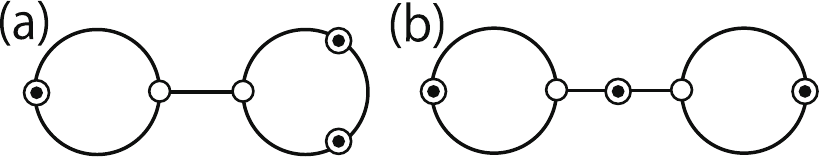}}
\caption{\label{Fig02}
 Graphs with three black vertices.}
\end{figure}

In this paper we shall show the following:

\begin{theorem}
\label{MainTheorem} 
Let $\Gamma$ be a minimal chart of type $(m;5,2)$.
Suppose that there exists a connected component of
$\Gamma_m$ with exactly five white vertices.
Then $\Gamma_m$ contains 
 one of the two graphs as shown in 
Fig.~\ref{Fig02}.
\end{theorem}

The paper is organized as follows.
In Section~\ref{s:Prel},
we define charts and minimal charts.
Let $\Gamma$ be a minimal chart, and $m$ a label of $\Gamma$. 
In Section~\ref{s:Lens},
we review a useful lemma for a disk called a lens.
In Section~\ref{s:kAngledDisk},
 we investigate a disk called 
a $k$-angled disk of $\Gamma_m$ with
at most one white vertex in its interior, 
where
a $k$-angled disk is a disk 
whose boundary contains exactly $k$ white vertices
and consists of edges of label $m$.
In Section~\ref{s:5AngledDisk},
we investigate a 5-angled disk of $\Gamma_m$.
In Section~\ref{s:4AngledDisk},
 we investigate a 4-angled disk of $\Gamma_m$.
In Section~\ref{s:TypeATypeC},
 we shall show that if $\Gamma$ is a minimal chart 
of type $(m; 5, 2)$,
then the graph $\Gamma_m$ contains neither graphs as shown 
in Fig.~\ref{Fig13}(a),(c).
In Section~\ref{s:TypeB},
 we shall show that if $\Gamma$ is a minimal chart 
of type $(m; 5, 2)$,
then the graph $\Gamma_m$ does not contain
the graph as shown 
in Fig.~\ref{Fig13}(b).
In Section~\ref{s:TypeD},
 we shall show that if $\Gamma$ is a minimal chart 
of type $(m; 5, 2)$,
then the graph $\Gamma_m$ does not contain
the graph as shown 
in Fig.~\ref{Fig13}(d).
In Section~\ref{s:IOC},
we review IO-Calculation(a property of numbers of 
inward arcs of label $k$ 
and outward arcs of label $k$ in a closed domain $F$
with $\partial F\subset\Gamma_{k-1}\cup\Gamma_k\cup\Gamma_{k+1}$
for some label $k$).
In Section~\ref{s:TypeE},
 we shall show that if $\Gamma$ is a minimal chart 
of type $(m; 5, 2)$,
then the graph $\Gamma_m$ does not contain
the graph as shown 
in Fig.~\ref{Fig13}(e).
In Section~\ref{s:TypeF},
 we shall show that if $\Gamma$ is a minimal chart 
of type $(m; 5, 2)$,
then the graph $\Gamma_m$ does not contain
the graph as shown 
in Fig.~\ref{Fig13}(f).
In Section~\ref{s:TriangleLemma},
we review Triangle Lemma. These lemmas will be used
in Section~\ref{s:TypeI}.
In Section~\ref{s:TypeI},
 we shall show that if $\Gamma$ is a minimal chart 
of type $(m; 5, 2)$,
then the graph $\Gamma_m$ does not contain
the graph as shown 
in Fig.~\ref{Fig13}(g).
Moreover, we shall prove Theorem~\ref{MainTheorem}.


\section{Preliminaries}
\label{s:Prel}

In this section, 
we introduce 
the definition of charts and its related words.

\begin{figure}[htb]
\begin{center}
\includegraphics{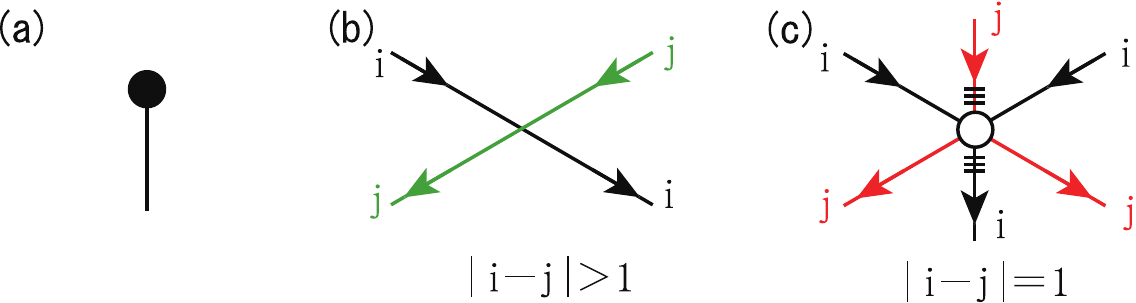}
\end{center}
\caption{ \label{Fig03} (a) A black vertex. (b) A crossing. (c) A white vertex. 
Each arc with three transversal short arcs is a middle arc at the white vertex. }
\end{figure}

Let $n$ be a positive integer.
An $n$-{\it chart}  
(a braid chart of degree $n$ \cite{KnottedSurfaces}
or a surface braid chart of degree $n$ \cite{BraidBook}) 
is 
an oriented labeled graph in the interior of a disk,
which may be empty 
or
have closed edges without vertices
satisfying the following four conditions
(see Fig.~\ref{Fig03}):
\begin{enumerate}
\item[(i)] 
Every vertex has degree $1$, $4$, or $6$.
\item[(ii)] 
The labels of edges are 
in $\{1,2,\dots,n-1\}$.
\item[(iii)]
In a small neighborhood of
each vertex of degree $6$,
there are six short arcs,
three consecutive arcs are
oriented inward 
and
the other three are outward,
and
these six are labeled $i$ and $i+1$
alternately for some $i$,
where the orientation and label of
each arc are inherited from
the edge containing the arc.
\item[(iv)]
For each vertex of degree $4$,
diagonal edges have the same label
and
are oriented coherently,
and the labels $i$ and $j$ of
the diagonals satisfy $|i-j|>1$.
\end{enumerate}
We call a vertex of degree $1$ a {\it black vertex},
a vertex of degree $4$ a {\it crossing}, and 
a vertex of degree $6$ a {\it white vertex}
respectively.

Among six short arcs
in a small neighborhood of
a white vertex,
a central arc of each three consecutive arcs
oriented inward (resp. outward) 
is called a   
{\it middle arc} at the white vertex
(see Fig.~\ref{Fig03}(c)).
For each white vertex $v$, 
there are two middle arcs at $v$ 
in a small neighborhood of $v$.
An edge is said to be {\it middle at} a white vertex $v$ if it contains a middle arc at $v$.

Let $e$ be an edge connecting $v_1$ and $v_2$.
If $e$ is oriented from $v_1$ to $v_2$,
then we say that 
$e$ is oriented {\it outward at $v_1$}
and {\it inward at $v_2$}.


\begin{figure}
\begin{center}
\includegraphics{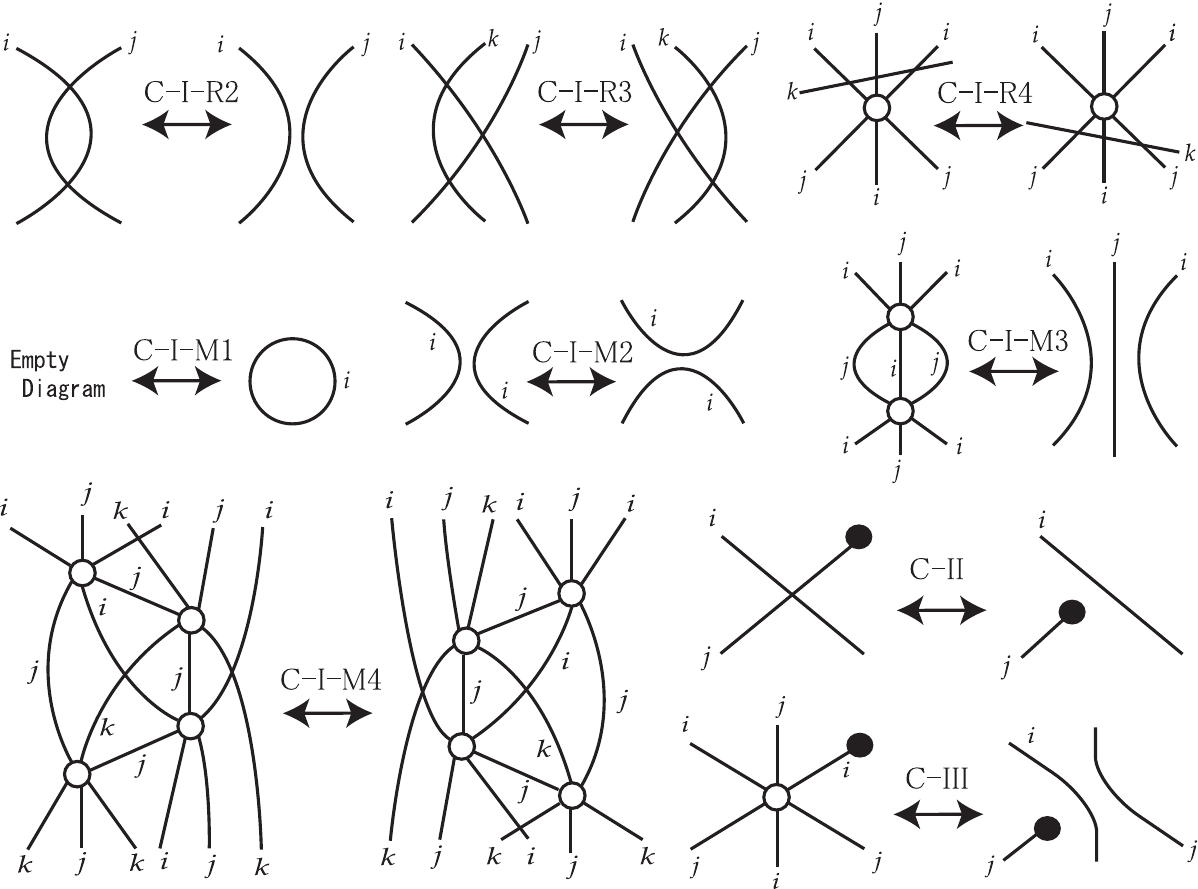}
\end{center}
\caption{ \label{Fig04} For the C-III move, 
the edge with the black vertex is not middle at
a white vertex in the left figure. }
\end{figure}

Now {\it C-moves} are local modifications 
of charts as shown in Fig.~\ref{Fig04}
(cf. \cite{KnottedSurfaces}, 
\cite{BraidBook} and \cite{Tanaka}).
Two charts are said to be {\it C-move equivalent}  if there exists
a finite sequence of C-moves 
which modifies one of the two charts 
to the other.

An edge in a chart is called 
a {\it free edge}
if it has
two black vertices.

For each chart $\Gamma$,
let $w(\Gamma)$ and $f(\Gamma)$ be the number of white vertices, and the number of free edges respectively.
The pair $(w(\Gamma), -f(\Gamma))$ is called a {\it complexity} of the chart (see \cite{BraidThree}).
A chart $\Gamma$ is called a {\it minimal chart} if its complexity is minimal among the charts C-move equivalent to the chart $\Gamma$ with respect to the lexicographic order of pairs of integers.

We showed the difference of a chart in a disk and in a 2-sphere (see \cite[Lemma 2.1]{ChartApp1}).
This lemma follows from that there exists a natural one-to-one correspondence between $\{$charts in $S^2\}/$C-moves and $\{$charts in $D^2\}/$C-moves, conjugations
(\cite[Chapter 23 and Chapter 25]{BraidBook}).
To make the argument simple, we assume that 
the charts lie on the 2-sphere instead of the disk.
\begin{assumption}
In this paper,
all charts are contained in the $2$-sphere $S^2$.
\end{assumption}
We have the special point in the 2-sphere $S^2$, called the point at infinity,
 denoted by $\infty$.
In this paper, all charts are contained in a disk such that the disk 
does not contain the point at infinity $\infty$.

Let $\Gamma$ be a chart,
and $m$ a label of $\Gamma$. 
A {\it hoop} is a closed edge of $\Gamma$ without vertices 
(hence without crossings, neither).
A {\it ring} is a simple closed curve in $\Gamma_m$ containing at least one crossing but not containing any white vertices.
A hoop is said to be {\it simple} 
if one of the two complementary domains
of the hoop
does not contain any white vertices.

We can assume that
all minimal charts $\Gamma$
satisfy the following four conditions 
(see \cite{ChartApp1},\cite{ChartAppII},\cite{ChartAppIII},\cite{StI}):

\begin{assumption}
\label{AssumeTerminal}
If an edge of $\Gamma$
contains a black vertex,
then the edge is a free edge 
or a terminal edge.
Moreover 
any terminal edge contains a middle arc.
\end{assumption}

\begin{assumption}
\label{NoSimpleHoop}
All free edges and simple hoops in $\Gamma$ 
are moved into a small neighborhood $U_\infty$ 
of the point at infinity $\infty$. 
Hence
we assume that 
$\Gamma$ does not contain free edges
nor simple hoops, 
otherwise mentioned. 
\end{assumption}

\begin{assumption}
\label{Ring}
Each complementary domain of
any ring and hoop must contain 
at least one white vertex. 
\end{assumption}

\begin{assumption}
\label{Infinity}
The point at infinity $\infty$ is moved in any complementary domain of $\Gamma$.
\end{assumption}

In this paper
for a subset $X$ in a space
we denote 
the interior of $X$,
the boundary of $X$ and
the closure of $X$
by Int$X$, $\partial X$
and $Cl(X)$
respectively.

Let $\alpha$ be a simple arc or an edge,
and $p,q$ the endpoints of $\alpha$.
We denote 
$\partial \alpha=\{p,q\}$
and ${\rm Int}\alpha=\alpha-\{p,q\}$.


\newpage

\section{Lenses}
\label{s:Lens}

In this section,
we review a useful lemma for a disk called a lens.

Let $\Gamma$ be a chart, 
and $m$ a label of $\Gamma$.
Let $L$ be the closure of a connected component 
of the set obtained by taking out 
all the white vertices from $\Gamma_m$.
If $L$ contains at least one white vertex
but does not contain any black vertex,
then $L$ is called an {\it internal edge of label $m$}.
Note that an internal edge may contain a crossing of $\Gamma$.

\begin{figure}[htb]
\centerline{\includegraphics{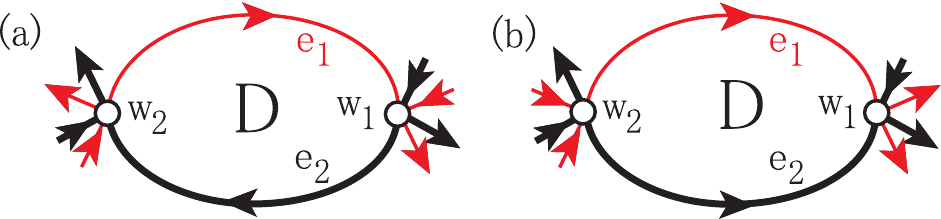}}
\caption{\label{Fig05}
Lenses.}
\end{figure}

Let $\Gamma$ be a chart. 
Let $D$ be a disk 
such that 
\begin{enumerate}
\item[(1)] the boundary $\partial D$ consists of an internal edge $e_1$ of label $m$ and an internal edge $e_2$ of label ${m+1}$, and 
\item[(2)] any edge containing a white vertex in $e_1$ does not intersect the open disk Int$D$.
\end{enumerate}
Note that $\partial D$ may contain crossings.
Let $w_1$ and $w_2$ be the white vertices in $e_1$. 
If the disk $D$ satisfies one of the following conditions, then $D$ is called  {\it a lens of type $(m,m+1)$}
(see Fig.~\ref{Fig05}):
\begin{enumerate}
	\item[(i)] Neither $e_1$ nor $e_2$ contains a middle arc. 
	\item[(ii)] One of the two edges $e_1$ and $e_2$ contains middle arcs at both white vertices $w_1$ and $w_2$ simultaneously.
\end{enumerate}

\begin{lemma}{\rm (\cite[Theorem 1.1]{ChartApp1})}
\label{LensThreeWhiteVertex}
 There exist at least three white vertices
in the interior of the lens for any minimal chart.
\end{lemma}

\begin{lemma}{\rm (\cite[Corollary 1.3]{ChartAppII})}
\label{NoLens}
 There is no lens in any minimal chart with 
at most seven white vertices.
\end{lemma}

Let $\Gamma$ be a chart,
and $m$ a label of $\Gamma$. 
A {\it loop} is a simple closed curve in $\Gamma_m$ with exactly one white vertex
(possibly with crossings).

\begin{figure}[thb]
\centerline{\includegraphics{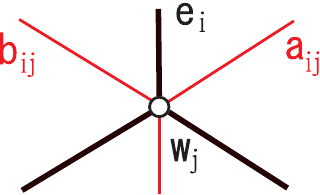}}
\caption{\label{Fig06}
The three edges $a_{ij},e_i,b_{ij}$ are consecutive edges around the white vertex $w_j$.}
\end{figure}

In our argument,
we often need a name for an unnamed edge by using a given edge and a given white vertex.
For the convenience,
we use the following naming:
Let $e',e_i,e''$ be three consecutive edges containing  a white vertex $w_j$. Here, 
the two edges $e'$ and $e''$ are unnamed edges. 
There are six arcs in a neighborhood $U$ of the white vertex $w_j$. 
If the three arcs $e'\cap U$, $e_i \cap U$, $e'' \cap U$ lie anticlockwise around the white vertex $w_j$ in this order, 
then $e'$ and $e''$ are denoted by $a_{ij}$ and $b_{ij}$ 
respectively (see Fig.~\ref{Fig06}).
There is a possibility $a_{ij}=b_{ij}$ if they are contained in a loop.

\begin{figure}
\begin{center}
\includegraphics{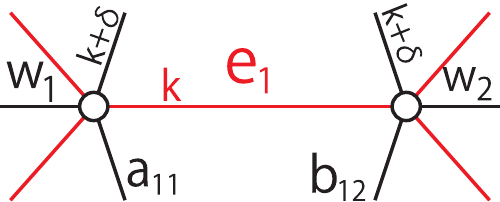}
\end{center}
\caption{\label{Fig07} The edge $e_1$ is of label $k$,
and $\delta\in\{+1,-1\}.$ }
\end{figure}

\begin{lemma}
\label{C-I-M2A11B12}
Let $\Gamma$ be a chart, and
$k$ a label of $\Gamma$.
Let $e_1$ be an internal edge of label $k$
with two white vertices $w_1$ and $w_2$
$($see Fig.~\ref{Fig07}$)$.
Suppose that $w_1,w_2\in\Gamma_{k+\delta}$
for some $\delta\in\{+1,-1\}$,
and suppose that one of the two edges $a_{11},b_{12}$
is a terminal edge.
If $\Gamma_{k+2\delta}\cap e_1=\emptyset$,
and if $\Gamma$ satisfies one of the following four conditions,
then $\Gamma$ is not a minimal chart. 
\begin{enumerate}
\item[{\rm (a)}]
The two edges $a_{11},b_{12}$ are oriented inward $($or outward$)$ at $w_1,w_2$, respectively.
\item[{\rm (b)}] 
The edge $a_{11}$ $($resp. $b_{12})$ is a terminal edge, and
$b_{12}$ $($resp. $a_{11})$ is not middle at the white vertex
different from $w_2$ $($resp. $w_{1})$.
\item[{\rm (c)}] 
The two edges $a_{11},b_{12}$ are middle at  $w_1,w_2$, respectively.
\item[{\rm (d)}] 
Both of  $a_{11},b_{12}$ are terminal edges.
\end{enumerate}
\end{lemma}

\begin{figure}
\begin{center}
\includegraphics{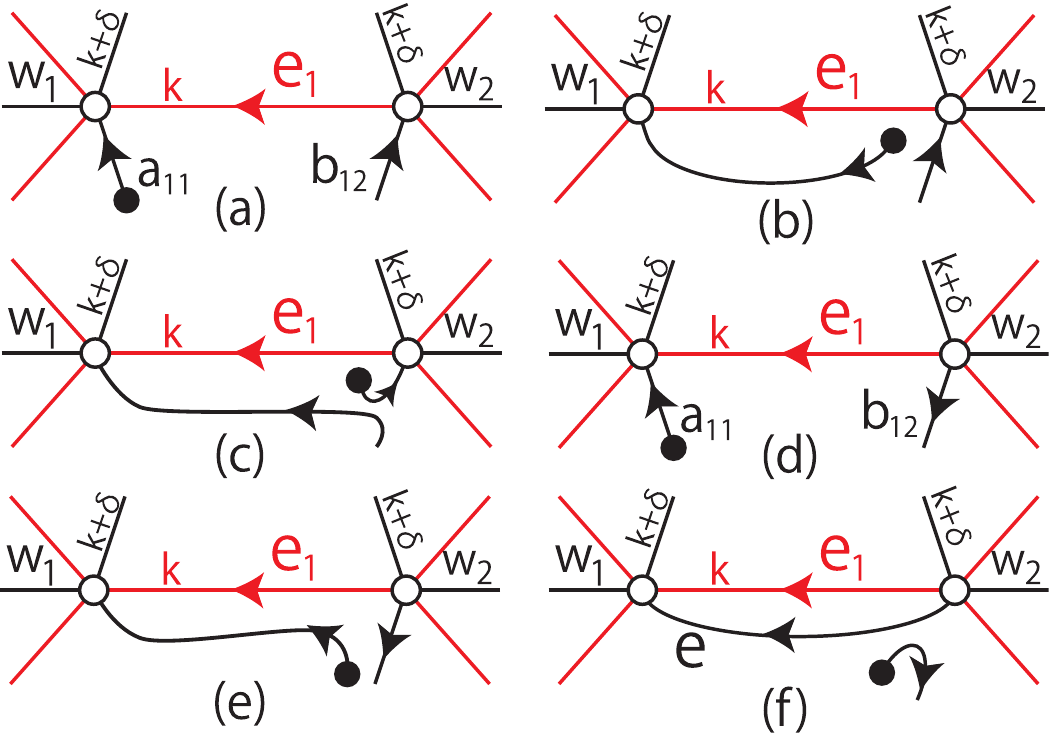}
\end{center}
\caption{\label{Fig08} The edge $e_1$ is of label $k$,
and $\delta\in\{+1,-1\}.$ }
\end{figure}

\begin{Proof}
Suppose that $\Gamma$ is a minimal chart.
Without loss of generality we can assume that
\begin{enumerate}
\item[(1)] 
$a_{11}$ is a terminal edge and oriented inward at $w_1$.
\end{enumerate}
Then by Assumption~\ref{AssumeTerminal},
the terminal edge $a_{11}$ is middle at $w_1$.
Thus 
\begin{enumerate}
\item[(2)] 
the edge $e_1$ is oriented inward at $w_1$ (i.e. 
the edge $e_1$ is oriented from $w_2$ to $w_1$).
\end{enumerate}

Since $e_1$ is an edge of label $k$,
we have
\begin{enumerate}
\item[(3)] $\Gamma_{k+\delta}\cap {\rm Int}e_1=\emptyset$. 
\end{enumerate}

Now by the condition of this lemma,
we have 
\begin{enumerate}
\item[(4)] $\Gamma_{k+2\delta}\cap e_1=\emptyset$. 
\end{enumerate}

First, we shall show that
if $\Gamma$ satisfies Condition~(a),
then we have a contradiction.
Since $a_{11}$ is oriented inward at $w_1$ by (1),
the edge $b_{12}$ is also oriented inward at $w_2$
by Condition~(a)
(see Fig.~\ref{Fig08}(a)).
Hence by (2)
\begin{enumerate}
\item[(5)] the edge $b_{12}$ is not middle at $w_2$.
\end{enumerate}

By (3),(4) and applying C-II moves along the edge $e_1$,
we can move the black vertex in $a_{11}$ 
near the white vertex $w_2$
 (see Fig.~\ref{Fig08}(b)).
Apply a C-I-M2 move between $a_{11}$ and $b_{12}$,
we obtain a new terminal edge of label $k+\delta$ at $w_2$
(see Fig.~\ref{Fig08}(c)).
However by (5), the terminal edge is not middle at $w_2$.
This contradicts Assumption~\ref{AssumeTerminal}.

Now by (1) and Lemma~\ref{C-I-M2A11B12}(a), we can assume that
\begin{enumerate}
\item[(6)] the edge $b_{12}$ is oriented outward at $w_2$
(see Fig.~\ref{Fig08}(d)).
\end{enumerate}

Next, we shall show that
if $\Gamma$ satisfies Condition~(b),
then we have a contradiction.
Let $w_3$ be the white vertex in $b_{12}$
different from $w_2$. Then by Condition~(b),
\begin{enumerate}
\item[(7)] the edge $b_{12}$ is not middle at $w_3$.
\end{enumerate}
By the similar way of the proof of Lemma~\ref{C-I-M2A11B12}(a),
we can move the black vertex in $a_{11}$
near the white vertex $w_2$ by C-II moves
(see Fig.~\ref{Fig08}(e)).
Apply a C-I-M2 move between $a_{11}$ and $b_{12}$,
we obtain an internal edge of label $k+\delta$ with $w_1$ and $w_2$,
and a terminal edge of label $k+\delta$ at $w_3$
(see Fig.~\ref{Fig08}(f)).
However we have the same contradiction by (7).

Next,  we shall show that
if $\Gamma$ satisfies Condition~(c),
then we have a contradiction.
By the similar way of the proof of Lemma~\ref{C-I-M2A11B12}(b),
we obtain an internal edge of label $k+\delta$
with $w_1$ and $w_2$, say $e$
(see Fig.~\ref{Fig08}(f)).
Thus by Condition~(c),
the edge $e$ is middle at both white vertices $w_1,w_2$.
Hence $e_1\cup e$ bounds a lens whose interior
does not contain any white vertices.
This contradicts Lemma~\ref{LensThreeWhiteVertex}.

Finally,
if $\Gamma$ satisfies Condition~(d),
then by Assumption~\ref{AssumeTerminal}
the two terminal edges $a_{11}$ and $b_{12}$
are middle at $w_1,w_2$,
respectively.
Thus $\Gamma$ satisfies Condition~(c).
Hence 
we have the same contradiction.

Therefore we have a contradiction for all cases.
Thus $\Gamma$ is not a minimal chart.
\end{Proof}

\section{Special $k$-angled disks}
\label{s:kAngledDisk}

In this section, we investigate a disk called a $k$-angled disk.

Let $X$ be a set in a chart $\Gamma$.
Let
 $$w(X)=\text{the number of white vertices in $X$.}$$

\begin{lemma}
\label{LemmaWithTerminal3}
{\rm (\cite[Lemma 3.2(1)]{ChartAppV})}
Let $\Gamma$ be a minimal chart,
and $m$ a label of $\Gamma$.
Let $G$ be a connected component of $\Gamma_m$.
If $1\le w(G)$, then $2\le w(G)$.
\end{lemma}

\begin{lemma}
$(${\rm \cite[Lemma 6.1]{ChartAppIII}}$)$
\label{ConditionRing} 
Let $\Gamma$ be a minimal chart.
Let $C$ be a ring or a non simple hoop, 
and $D$ a disk with $\partial D=C$.
If $w(\Gamma\cap D)=1$, then $\Gamma$ is C-move equivalent to the minimal chart $Cl(\Gamma-C)$. 
\end{lemma}

Let $\Gamma$ be a chart, $m$ a label of $\Gamma$, 
$D$ a disk with $\partial D\subset \Gamma_m$, 
and $k$ a positive integer.
If $\partial D$ contains exactly
$k$ white vertices, 
then $D$ is called 
{\it a $k$-angled disk of $\Gamma_m$}. 
Note that 
the boundary $\partial D$ may contain crossings.

\begin{lemma}
\label{GammaM-1MissD}
Let $\Gamma$ be a minimal chart, and
$m$ a label of $\Gamma$.
Let $D$ be a $k$-angled disk of $\Gamma_m$.
Suppose that all the white vertices on $\partial D$
are contained in $\Gamma_{m+\varepsilon}$ 
for some $\varepsilon\in\{+1,-1\}$.
If $w(\Gamma\cap {\rm Int}D)\leqq1$,
then $\Gamma$ can be modified to a minimal chart $\Gamma'$
by C-moves in ${\rm Int}D$ such that
$\Gamma'_{m-\varepsilon}\cap D=\emptyset$.
\end{lemma}

\begin{Proof}
We shall show that ${\rm Int}D$ does not contain any white vertex
in $\Gamma_{m-\varepsilon}$.
Suppose that
there exists a white vertex $w$ in $\Gamma_{m-\varepsilon}$.
Let $G$ be the connected component of $\Gamma_{m-\varepsilon}$
with $G\ni w$.
Then by the condition of this lemma,
we have $G\subset {\rm Int}D$.
Thus the condition $w(\Gamma\cap {\rm Int}D)\leqq1$
implies $w(G)=1$.
This contradicts Lemma~\ref{LemmaWithTerminal3}.
Hence  ${\rm Int}D$ does not contain any white vertex
in $\Gamma_{m-\varepsilon}$.

Suppose that $\Gamma_{m-\varepsilon}\cap D\not=\emptyset$.
Then we shall show that the set $\Gamma_{m-\varepsilon}\cap D$
consists of rings and non simple hoops.
Let $e$ be an internal edge (possibly a ring or a hoop)
of label $m-\varepsilon$ intersecting the disk $D$.
Since all the white vertices on $\partial D$
are contained in  $\Gamma_{m+\varepsilon}$,
the condition $\partial D\subset \Gamma_m$
implies $e\subset {\rm Int}D$.
Since ${\rm Int}D$ does not contain any white vertex
in $\Gamma_{m-\varepsilon}$,
the edge $e$ is a ring or a hoop.
Thus by Assumption~\ref{Ring},
the curve $e$ is a ring or a non simple hoop.

Since $w(\Gamma\cap{\rm Int}D)\leqq 1$, 
by Lemma~\ref{ConditionRing}
we can eliminate each ring and non simple hoop in ${\rm Int}D$.
Hence the chart $\Gamma$ 
can be modified to a minimal chart $\Gamma'$
by C-moves in ${\rm Int}D$ such that
$\Gamma'_{m-\varepsilon}\cap D=\emptyset$.
\end{Proof}

Let $\Gamma$ be a chart, and
$m$ a label of $\Gamma$.
An edge of label $m$ is called a {\it feeler} of a $k$-angled disk $D$ of $\Gamma_m$
if the edge intersects $N-\partial D$
where $N$ is a regular neighborhood of $\partial D$ in $D$.

Let $\Gamma$ be a chart,
and $D$ a $k$-angled disk of $\Gamma_m$.
If any feeler of $D$ of label $m$ is a terminal edge,
then $D$ is called a {\it special} $k$-angled disk.

\begin{figure}
\centerline{\includegraphics{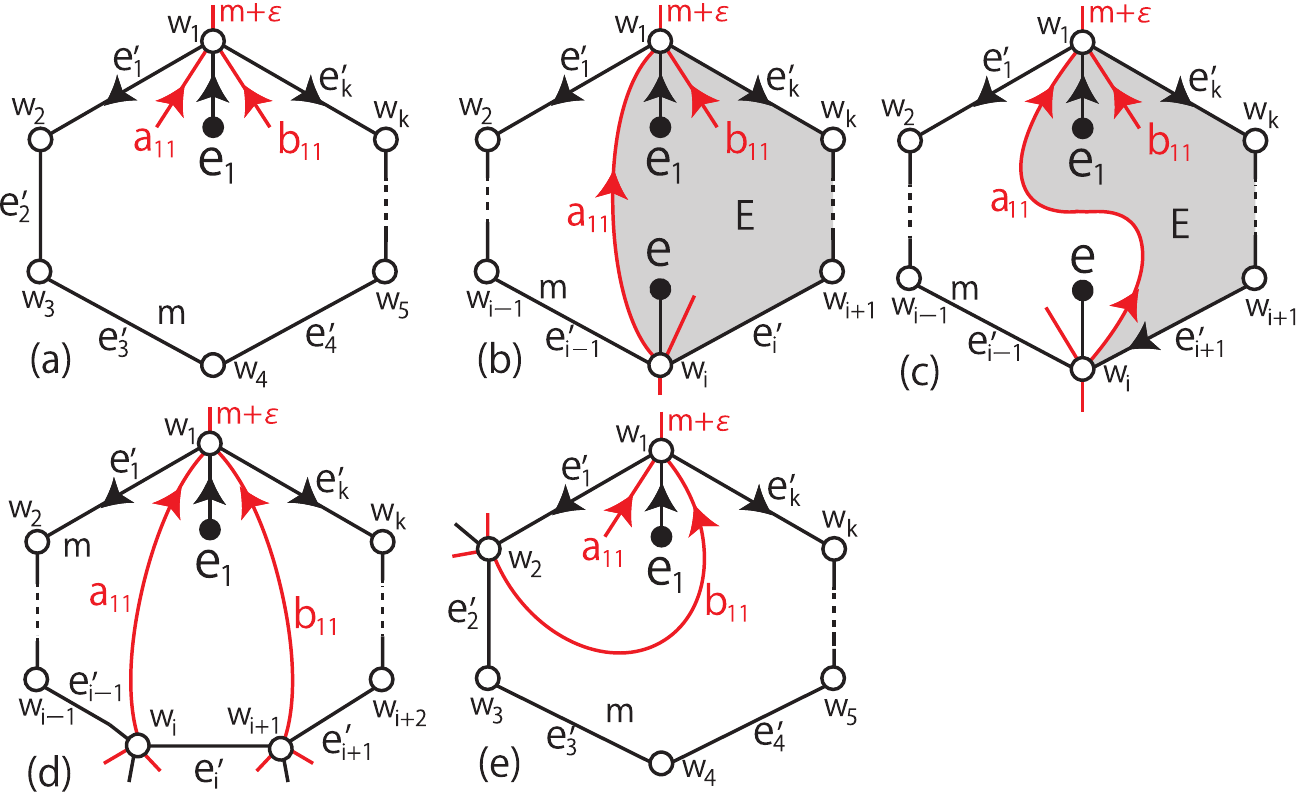}}
\caption{\label{Fig09} 
The edge $e_1$ is a terminal edge of label $m$,
and $a_{11},b_{11}$ are internal edges of label $m+\varepsilon$
for some $\{+1,-1\}$. 
The gray regions are disks $E$.}
\end{figure}

\begin{lemma}
\label{kAngledLessOneWhite}
Let $\Gamma$ be a minimal chart, and 
$m$ a label of $\Gamma$.
Let $D$ be a special $k$-angled disk of $\Gamma_m$
with at least one feeler $e_1$.
Let $w_1$ be the white vertex in $e_1$,
and $w_2,w_3,\cdots,w_k$ the other white vertices on $\partial D$
lying anticlockwise in this order.
Suppose that all of $w_1,w_2,\cdots,w_k$
are contained in $\Gamma_{m+\varepsilon}$
for some $\varepsilon\in\{+1,-1\}$.
Let $a_{11},b_{11}$ be internal edges 
$($possibly terminal edges$)$ of label $m+\varepsilon$ 
at $w_1$ in $D$ such that
$a_{11},e_1,b_{11}$ lie anticlockwise in this order
$($see Fig.~\ref{Fig09}$($a$))$.
If $w(\Gamma \cap{\rm Int}D)\leqq1$,
then the following five conditions hold:

\begin{enumerate}
\item[{\rm (a)}] 
If $a_{11}\ni w_i$ or $b_{11}\ni w_i$ for some $i\in\{2,3,\cdots,k\}$,
then $D$ does not contain a feeler at $w_i$.
\item[{\rm (b)}]
$a_{11}\not\ni w_2$ and $b_{11}\not\ni w_k$. 
\item[{\rm (c)}] 
If $a_{11}\ni w_i$ for some $i\in\{3,4,\cdots,k-1\}$,
then $b_{11}\not\ni w_{i+1}$.
\item[{\rm (d)}] 
If $b_{11}\ni w_i$ for some $i\in\{3,4,\cdots,k-1\}$,
then $a_{11}\not\ni w_{i-1}$.
\item[{\rm (e)}] 
$b_{11}\not\ni w_2$ and $a_{11}\not\ni w_k$. 
\end{enumerate}
\end{lemma}

\begin{Proof}
Let $e_1',e_2',\cdots,e_k'$ be the internal edges
of label $m$ in $\partial D$
such that 
$\partial e_i'=\{w_i,w_{i+1}\}$ for $i=1,2,\cdots,k-1$
and
$\partial e_k'=\{w_k,w_1\}$.
Without loss of generality, we can assume that
\begin{enumerate}
\item[(1)]
the terminal edge $e_1$ is oriented inward at $w_1$.
\end{enumerate}
Then by Assumption~\ref{AssumeTerminal},
\begin{enumerate}
\item[(2)]
both of $a_{11},b_{11}$ are oriented inward at $w_1$ and 
\item[(3)]
both of $e_1',e_k'$ are oriented outward at $w_1$
(see Fig.~\ref{Fig09}(a)). 
\end{enumerate}

Since $w(\Gamma \cap{\rm Int}D)\leqq1$, 
by Lemma~\ref{GammaM-1MissD}
we can assume that 
\begin{enumerate}
\item[(4)]
$\Gamma_{m-\varepsilon}\cap D=\emptyset$
($\Gamma_{m-\varepsilon}\cap a_{11}=\emptyset$
and $\Gamma_{m-\varepsilon}\cap b_{11}=\emptyset$).
\end{enumerate}

First, we shall show Statement~(a).
Now, suppose that
$a_{11}\ni w_i$ for some $i\in\{2,3,\cdots,k\}$.
Then the edge $a_{11}$ separates the disk $D$
into two disks.
One of the two disks contains the terminal edge $e_1$,
say $E$.
Suppose that $D$ contains a feeler $e$ at $w_i$.
Then $e\subset E$ or $e\not\subset E$.

If $e\subset E$ (see Fig.~\ref{Fig09}(b)), 
then by (4) and Lemma~\ref{C-I-M2A11B12}(d)
the chart $\Gamma$ is not minimal. 
This contradicts the fact that $\Gamma$ is minimal.
Thus $e\not\subset E$ 
(see Fig.~\ref{Fig09}(c)).

By (2),
the edge $a_{11}$ is oriented outward at $w_i$.
Thus
by Assumption~\ref{AssumeTerminal},
 the edge $e_i'$ is oriented inward at $w_i$.
Hence by (1),(4) and Lemma~\ref{C-I-M2A11B12}(a)
the chart $\Gamma$ is not minimal. 
This contradicts the fact that $\Gamma$ is minimal.
Therefore $D$ does not contain a feeler at $w_i$.

Similarly if $b_{11}\ni w_i$,
then we can show  that 
$D$ does not contain a feeler at $w_i$.
Thus Statement~(a) holds.

Next, we shall show Statement~(b).
If $a_{11}\ni w_2$,
then
by Lemma~\ref{kAngledLessOneWhite}(a)
the disk $D$ does not contain a feeler at $w_2$.
Thus the curve $e_1'\cup a_{11}$ bounds a lens whose interior
contains at most one white vertex.
This contradicts Lemma~\ref{LensThreeWhiteVertex}.
Thus $a_{11}\not\ni w_2$.

Similarly we can show $b_{11}\not\ni w_k$.
Thus Statement~(b) holds.

Next, we shall show Statement~(c).
Now, suppose that 
$a_{11}\ni w_i$ for some $i\in\{3,4,\cdots,k-1\}$.
We shall show $b_{11}\not\ni w_{i+1}$. 
If $b_{11}\ni w_{i+1}$,
then by Lemma~\ref{kAngledLessOneWhite}(a),
the disk $D$ contains neither a feeler at $w_i$ nor
a feeler at $w_{i+1}$
(see Fig.~\ref{Fig09}(d)). 

Now, by (1), the terminal edge $e_1$ is oriented inward at $w_1$.
Since the edge $e_i'$ is oriented inward at $w_i$ or $w_{i+1}$,
the chart $\Gamma$ is not minimal 
by (4) and Lemma~\ref{C-I-M2A11B12}(a).
This contradicts the fact that $\Gamma$ is minimal.
Therefore $b_{11}\not\ni w_{i+1}$. 
Thus Statement~(c) holds.

Similarly we can show Statement~(d).

Finally, we shall show Statement~(e).
Suppose $b_{11}\ni w_2$.
Then by Lemma~\ref{kAngledLessOneWhite}(a),
the disk $D$ does not contain a feeler at $w_2$
(see Fig.~\ref{Fig09}(e)).
Since the terminal edge $e_1$ is oriented inward at $w_1$
by (1),
and since the edge $e_1'$ is oriented inward at $w_2$ by (3),
 the chart $\Gamma$ is not minimal
by (4) and Lemma~\ref{C-I-M2A11B12}(a). 
This contradicts the fact that $\Gamma$ is minimal.
Therefore $b_{11}\not\ni w_2$.

Similarly we can show $a_{11}\not\ni w_k$.
Thus Statement~(e) holds.
\end{Proof}



\section{Special 5-angled disks}
\label{s:5AngledDisk}

In this section, we investigate a special $5$-angled disk whose 
interior contains at most one white vertex.

Let $\Gamma$ be a chart. 
Suppose that an object consists of 
some edges of $\Gamma$, arcs in edges of 
$\Gamma$ and arcs around white vertices.
Then the object is called a {\it pseudo chart}.

\begin{figure}
\centerline{\includegraphics{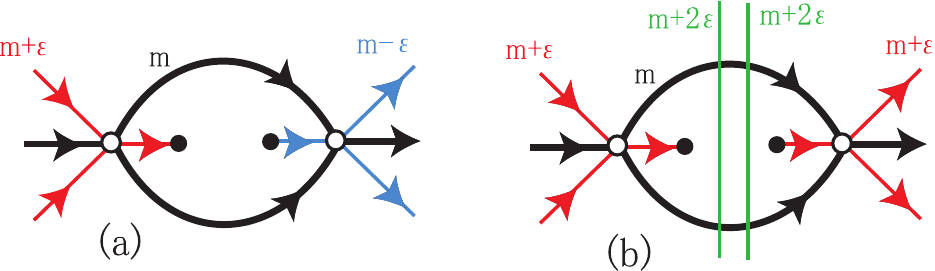}}
\caption{\label{Fig10}
$m$ is a label,
and $\varepsilon\in\{+1,-1\}$.}
\end{figure}

\begin{lemma}
\label{Theorem2AngledDisk}
{\rm (\cite[Corollary 5.8]{ChartAppII})}
Let $\Gamma$ be a minimal chart.
Let $D$ be a $2$-angled disk of $\Gamma_m$ with at most one feeler.
If $w(\Gamma\cap{\rm Int}D)=0$,
then a regular neighborhood of $D$ contains one of two pseudo charts as shown in Fig.~\ref{Fig10}.
\end{lemma}

\begin{lemma}
\label{5AngledOneWhite}
Let $\Gamma$ be a minimal chart, and 
$m$ a label of $\Gamma$.
Let $D$ be a special $5$-angled disk of $\Gamma_m$
with at least one feeler $e_1$.
Suppose that all the white vertices on $\partial D$
are contained in $\Gamma_{m+\varepsilon}$
for some $\varepsilon\in\{+1,-1\}$.
Let $w_1$ be the white vertex in $e_1$.
Let $a_{11},b_{11}$ be internal edges 
$($possibly terminal edges$)$ of label $m+\varepsilon$
at $w_1$ in $D$.
If $w(\Gamma \cap{\rm Int}D)\leqq 1$,
then one of $a_{11},b_{11}$ contains a white vertex 
in ${\rm Int}D$.
\end{lemma}

\begin{Proof}
Let $w_2,\cdots,w_5$ be the four white vertices on $\partial D$
different from $w_1$ such that $w_1,w_2,\cdots,w_5$ 
lie anticlockwise in this order.
Without loss of generality
we can assume that
$a_{11},e_1,b_{11}$ lie anticlockwise in this order.
By Assumption~\ref{AssumeTerminal},
\begin{enumerate}
\item[(1)]
neither $a_{11}$ nor $b_{11}$ is a terminal edge.
\end{enumerate}

Suppose that neither $a_{11}$ nor $b_{11}$
contains a white vertex in ${\rm Int}D$.
Then by applying Lemma~\ref{kAngledLessOneWhite}(b),(e) 
for the edge $a_{11}$, we have 
$a_{11}\ni w_3$ or $a_{11}\ni w_4$.

If $a_{11}\ni w_3$, then by (1)
we have $b_{11}\ni w_4$ or $b_{11}\ni w_5$.
This contradicts Lemma~\ref{kAngledLessOneWhite}(b),(c).

If $a_{11}\ni w_4$, then by (1)
we have $b_{11}\ni w_5$.
This contradicts Lemma~\ref{kAngledLessOneWhite}(b).

Therefore we have a contradiction for both cases.
Hence one of $a_{11},b_{11}$ contains a white vertex 
in ${\rm Int}D$.
\end{Proof}

Let $\Gamma$ be a chart, 
$D$ a disk, and 
$G$ a pseudo chart with $G \subset D$.
Let $r:D\to D$ be a reflection of $D$, and $G^*$ the pseudo chart obtained from $G$ by changing the orientations of all of the edges.
Then the set $\{G,G^*, r(G), r(G^*)\}$ 
is called the {\it RO-family of the pseudo chart $G$}.

\begin{figure}
\centerline{\includegraphics{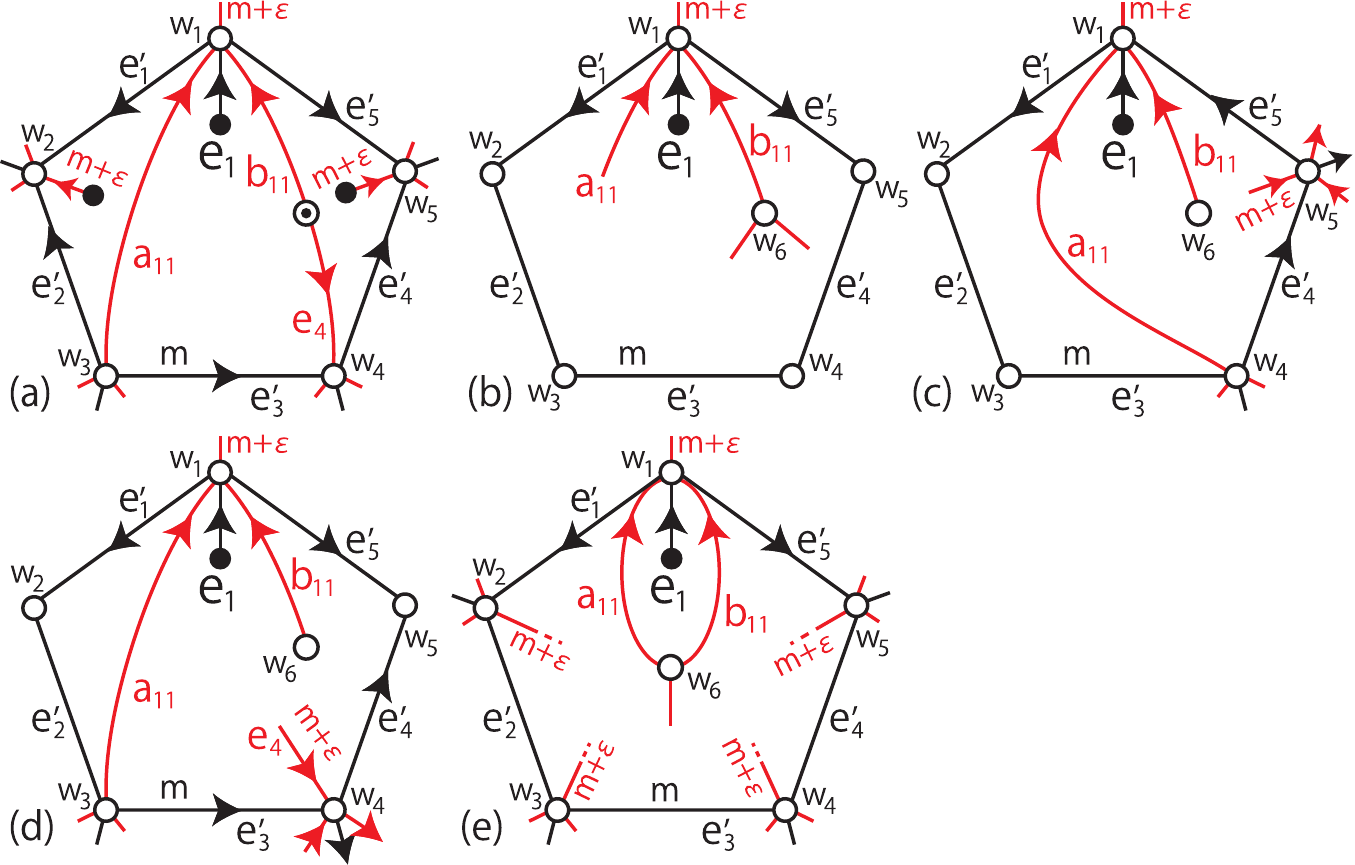}}
\caption{\label{Fig11} 
Special 5-angled disks with one feeler $e_1$.}
\end{figure}

\begin{lemma}
\label{5AngledOneFeelerOneWhiteCase1}
Let $\Gamma$ be a minimal chart, and 
$m$ a label of $\Gamma$.
Let $D$ be a special $5$-angled disk of $\Gamma_m$
with at least one feeler $e_1$.
Suppose that all the white vertices on $\partial D$
are contained in $\Gamma_{m+\varepsilon}$
for some $\varepsilon\in\{+1,-1\}$.
Let $w_1$ be the white vertex in $e_1$.
Let $a_{11},b_{11}$ be internal edges 
$($possibly terminal edges$)$ of label $m+\varepsilon$
at $w_1$ in $D$.
Suppose that one of $a_{11},b_{11}$ contains 
a white vertex in ${\rm Int}D$, but the other
 contains a white vertex in $\partial D$
different from $w_1$.
If $w(\Gamma \cap{\rm Int}D)\leqq 1$,
then the disk $D$ contains exactly one feeler $e_1$.
Moreover, 
the disk $D$ contains one of the RO-family
of the  pseudo chart as shown in 
Fig.~\ref{Fig11}$($a$)$.
\end{lemma}

\begin{Proof}
Let $w_2,\cdots,w_5$ be the four white vertices on $\partial D$
different from $w_1$ such that $w_1,w_2,\cdots,w_5$ 
lie anticlockwise in this order.
Let $e_1',e_2',\cdots,e_5'$ be the five internal edges of label $m$ in $\partial D$
with $\partial e_i'=\{w_i,w_{i+1}\}$ ($i=1,2,3,4$) and
$\partial e_5'=\{w_5,w_1 \}$.
Without loss of generality
we can assume that 
$a_{11},e_1,b_{11}$ lie anticlockwise in this order
(see Fig.~\ref{Fig11}(b)).

Since $w(\Gamma \cap{\rm Int}D)\leqq1$, 
by Lemma~\ref{GammaM-1MissD}
we can assume that 
\begin{enumerate}
\item[(1)]
$\Gamma_{m-\varepsilon}\cap D=\emptyset$
($\Gamma_{m-\varepsilon}\cap a_{11}=\emptyset$).
\end{enumerate}

Without loss of generality,
we can assume that
the terminal edge $e_1$ is oriented inward at $w_1$.
Then by Assumption~\ref{AssumeTerminal},
\begin{enumerate}
\item[(2)]
both of $a_{11},b_{11}$ are oriented inward at $w_1$,
\item[(3)] 
both of $e_1',e_5'$ are oriented outward at $w_1$.
\end{enumerate}
Without loss of generality,
we can assume that the edge $b_{11}$ contains a white vertex in ${\rm Int}D$, say $w_6$.
Thus by the condition of this lemma,
the edge $a_{11}$ contains a white vertex 
in $\partial D$
different from $w_1$.
Hence by 
Lemma~\ref{kAngledLessOneWhite}(b),(e), 
we have $a_{11}\ni w_3$ or $a_{11}\ni w_4$.

We shall show that $a_{11}\ni w_3$.
If $a_{11}\ni w_4$,
then by Lemma~\ref{kAngledLessOneWhite}(a)
the disk $D$ does not contain a feeler at $w_4$.
Hence by (1) and Lemma~\ref{C-I-M2A11B12}(a)
the internal edge $e_4'$ is oriented from $w_4$ to $w_5$
(because if not,
then the two edges $e_1,e_4'$ are oriented inward at $w_1,w_4$,
respectively.
Thus by (1) and Lemma~\ref{C-I-M2A11B12}(a),
the chart $\Gamma$ is not minimal.
This is a contradiction.).
Moreover, by (1) and Lemma~\ref{C-I-M2A11B12}(b),
the edge $e_4'$ is middle at $w_5$.
Hence by Assumption~\ref{AssumeTerminal},
the disk $D$ does not contain a feeler at $w_5$.
Thus by the definition of the chart,
the internal edge $e_5'$ is oriented from $w_5$ to $w_1$.
Hence the edge $e_5'$ must be oriented inward at $w_1$
(see Fig.~\ref{Fig11}(c)).
This contradicts (3).
Thus $a_{11}\ni w_3$.

Similarly,
by (1) and Lemma~\ref{C-I-M2A11B12}(a),(b),
we can show that 
the disk $D$ contains neither feeler at $w_3$ nor
feeler at $w_4$, and
\begin{enumerate}
\item[(4)]
the internal edge $e_3'$ is oriented from $w_3$ to $w_4$,
and is middle at $w_4$, 
\item[(5)]
the internal edge $e_4'$ is oriented from $w_4$ to $w_5$
(see Fig.~\ref{Fig11}(d)).
\end{enumerate}
Let $e_4$ be an internal edge 
(possibly a terminal edge) of label $m+\varepsilon$ 
at $w_4$ in $D$.
Then by (4),
\begin{enumerate}
\item[(6)] the edge $e_4$ is oriented inward at $w_4$,
but not middle at $w_4$.
\end{enumerate}
Thus by Assumption~\ref{AssumeTerminal},
we have $e_4\ni w_5$ or $e_4\ni w_6$.

We shall show that $e_4\ni w_6$.
If $e_4\ni w_5$,
then the disk $D$ does not contain a feeler at $w_5$
by Lemma~\ref{kAngledLessOneWhite}(b),(e).
Thus there are three consecitive edges $e_4',e_4,e_5'$ at $w_5$
such that $e_4',e_5'$ are oriented inward at $w_5$ 
by (3) and (5),
but $e_4$ is oriented outward at $w_5$ by (6).
This contradicts the definition of the chart.
Hence $e_4\ni w_6$.

Let $e_6$ be an internal edge (possibly a terminal edge)
of label $m+\varepsilon$ at $w_6$
different from $b_{11},e_4$.
Then by (2) and (6),
\begin{enumerate}
\item[(7)]
the edge $e_6$ is oriented inward at $w_6$.
\end{enumerate}

We shall show that
the disk $D$ does not contain a feeler at $w_5$.
If $D$ contains a feeler at $w_5$,
then two internal edges $e,e'$ of label $m+\varepsilon$ at $w_5$ in $D$
contain white vertices in ${\rm Int}D$
by Assumption~\ref{AssumeTerminal}.
Hence the condition $w({\rm Int}D)\leqq1$
implies that the two edges $e,e'$
contain the same white vertex $w_6$.
Thus 
there exist four internal edges $b_{11},e_4,e,e'$ 
of label $m+\varepsilon$
at $w_6$.
This contradicts the definition of the chart.
Hence $D$ does not contain a feeler at $w_5$.

Let $e_5$ be an internal edge (possibly a terminal edge)
of label $m+\varepsilon$ at $w_5$ in $D$.
Since both of $e_4',e_5'$ are oriented inward at $w_5$ by (3) and (5),
the edge $e_5$ is oriented inward at $w_5$.
Thus by (7), we have $e_5\not=e_6$.
Hence both of $e_5,e_6$ are terminal edges.

Moreover by Assumption~\ref{AssumeTerminal},
we can show that
the disk $D$ does not contain a feeler at $w_2$
and there exists a terminal edge of label $m+\varepsilon$ 
at $w_2$ in $D$.
Therefore
the disk $D$ contains 
the pseudo chart as shown in 
Fig.~\ref{Fig11}(a).
\end{Proof}


\begin{lemma}
\label{5AngledOneFeelerOneWhite}
Let $\Gamma$ be a minimal chart, and 
$m$ a label of $\Gamma$.
Let $D$ be a special $5$-angled disk of $\Gamma_m$
with at least one feeler.
Suppose that all the white vertices on $\partial D$
are contained in $\Gamma_{m+\varepsilon}$
for some $\varepsilon\in\{+1,-1\}$.
If  $w(\Gamma \cap{\rm Int}D)\leqq 1$,
then the disk $D$ contains exactly one feeler.
Moreover, 
the disk $D$ contains one of RO-families
of the two pseudo charts as shown in 
Fig.~\ref{Fig11}$($a$)$,$($e$)$.
\end{lemma}

\begin{Proof}
Let $e_1$ be a feeler of $D$,
and $w_1$ the white vertex in $e_1$.
Let $w_2,\cdots,w_5$ be the four white vertices on $\partial D$
different from $w_1$ such that $w_1,w_2,\cdots,w_5$ 
lie anticlockwise in this order.
Let $a_{11},b_{11}$ be internal edges (possibly terminal edges)
of label $m+\varepsilon$ at $w_1$ such that
$a_{11},e_1,b_{11}$ lie anticlockwise in this order
(see Fig.~\ref{Fig11}(b)).
By Lemma~\ref{5AngledOneWhite},
one of $a_{11},b_{11}$ contains a white vertex 
in ${\rm Int}D$.
Without loss of generality,
we can assume that 
the edge $b_{11}$ contains a white vertex in ${\rm Int}D$.
Moreover, by Assumption~\ref{AssumeTerminal},
the edge $a_{11}$ is not a terminal edge.
Hence either $a_{11}$ contains a white vertex in ${\rm Int}D$,
or $a_{11}$ contains a white vertex in $\partial D$
different from $w_1$.

If $a_{11}$ contains a white vertex in $\partial D$ 
different from $w_1$,
then by Lemma~\ref{5AngledOneFeelerOneWhiteCase1}
the disk $D$ contains one of the RO-family of
the pseudo chart as shown in Fig.~\ref{Fig11}(a).

Suppose that $a_{11}$ contains a white vertex in ${\rm Int}D$.
Since $w(\Gamma \cap{\rm Int}D)\leqq 1$,
both of $a_{11},b_{11}$ contain the same white vertex
in ${\rm Int}D$. 
Thus $a_{11}\cup b_{11}$ bounds a 2-angled disk in $D$.
Hence by Lemma~\ref{Theorem2AngledDisk},
the 2-angled disk has no feeler.

Finally we shall show that $D$ has exactly one feeler.
If $D$ has another feeler $e_2$ at some $w_i$ ($i=2,3,4,5$),
then $D$ has at least two feelers. 
Hence by Lemma~\ref{5AngledOneWhite} 
and Lemma~\ref{5AngledOneFeelerOneWhiteCase1},
the two internal edges $e,e'$ of label $m+\varepsilon$ 
at $w_i$ in $D$ contain white vertices in ${\rm Int}D$.
However, the condition $w(\Gamma \cap{\rm Int}D)\leqq 1$
implies that there are four edges $a_{11},b_{11},e,e'$ of 
label $m+\varepsilon$ at $w_6$.
This contradicts the definition of the chart.
Thus $D$ has exactly one feeler $e_1$.
Hence $D$ contains the pseudo chart as shown in 
Fig.~\ref{Fig11}(e).
\end{Proof}

From the above three lemmas,
we have the following corollary:

\begin{corollary}
\label{5AngledTwoFeeler}
Let $\Gamma$ be a minimal chart, and 
$m$ a label of $\Gamma$.
Let $D$ be a special $5$-angled disk of $\Gamma_m$.
Suppose that all the white vertices on $\partial D$
are contained in $\Gamma_{m+\varepsilon}$
for some $\varepsilon\in\{+1,-1\}$.
Then we have the following:
\begin{enumerate}
\item[{\rm (a)}]
If the disk $D$ contains at least one feeler,
then $w(\Gamma \cap{\rm Int}D)\geqq 1$.
\item[{\rm (b)}]
If the disk $D$ contains at least two feelers,
then $w(\Gamma \cap{\rm Int}D)\geqq 2$.
\end{enumerate}
\end{corollary}



\section{Special 4-angled disks}
\label{s:4AngledDisk}

In this section, we investigate a $4$-angled disk whose interior
contains at most one white vertex.
By the similar way of the proof of 
Lemma~\ref{5AngledOneWhite},
we can show the following lemma:

\begin{lemma}
\label{4AngledOneWhite}
Let $\Gamma$ be a minimal chart, and 
$m$ a label of $\Gamma$.
Let $D$ be a special $4$-angled disk of $\Gamma_m$
with at least one feeler $e_1$.
Suppose that all the white vertices on $\partial D$
are contained in $\Gamma_{m+\varepsilon}$
for some $\varepsilon\in\{+1,-1\}$.
Let $w_1$ be the white vertex in $e_1$.
Let $a_{11},b_{11}$ be internal edges 
$($possibly terminal edges$)$ of label $m+\varepsilon$
at $w_1$ in $D$.
If $w(\Gamma \cap{\rm Int}D)\leqq 1$,
then one of $a_{11},b_{11}$ contains a white vertex 
in ${\rm Int}D$.
\end{lemma}

\begin{figure}
\centerline{\includegraphics{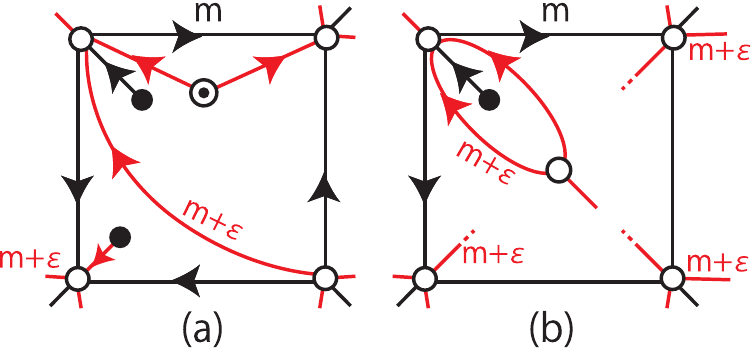}}
\caption{\label{Fig12} 
Special 4-angled disks with one feeler.}
\end{figure}

By the similar way of the proof of
Lemma~\ref{5AngledOneFeelerOneWhite},
we can show the following lemma:

\begin{lemma}
\label{4AngledOneFeelerOneWhite}
Let $\Gamma$ be a minimal chart, and 
$m$ a label of $\Gamma$.
Let $D$ be a special $4$-angled disk of $\Gamma_m$
with at least one feeler.
Suppose that all the white vertices on $\partial D$
are contained in $\Gamma_{m+\varepsilon}$
for some $\varepsilon\in\{+1,-1\}$.
If  $w(\Gamma \cap{\rm Int}D)\leqq 1$,
then the disk $D$ contains exactly one feeler.
Moreover, 
the disk $D$ contains one of RO-families
of the two pseudo charts as shown in 
Fig.~\ref{Fig12}.
\end{lemma}

From the above two lemmas,
we have the following corollary:

\begin{corollary}
\label{4AngledTwoFeeler}
Let $\Gamma$ be a minimal chart, and 
$m$ a label of $\Gamma$.
Let $D$ be a special $4$-angled disk of $\Gamma_m$.
Suppose that all the white vertices on $\partial D$
are contained in $\Gamma_{m+\varepsilon}$
for some $\varepsilon\in\{+1,-1\}$.
Then we have the following:
\begin{enumerate}
\item[{\rm (a)}]
If the disk $D$ contains at least one feeler,
then $w(\Gamma \cap{\rm Int}D)\geqq 1$.
\item[{\rm (b)}]
If the disk $D$ contains at least two feelers,
then $w(\Gamma \cap{\rm Int}D)\geqq 2$.
\end{enumerate}
\end{corollary}



\begin{figure}[htb]
\centerline{\includegraphics{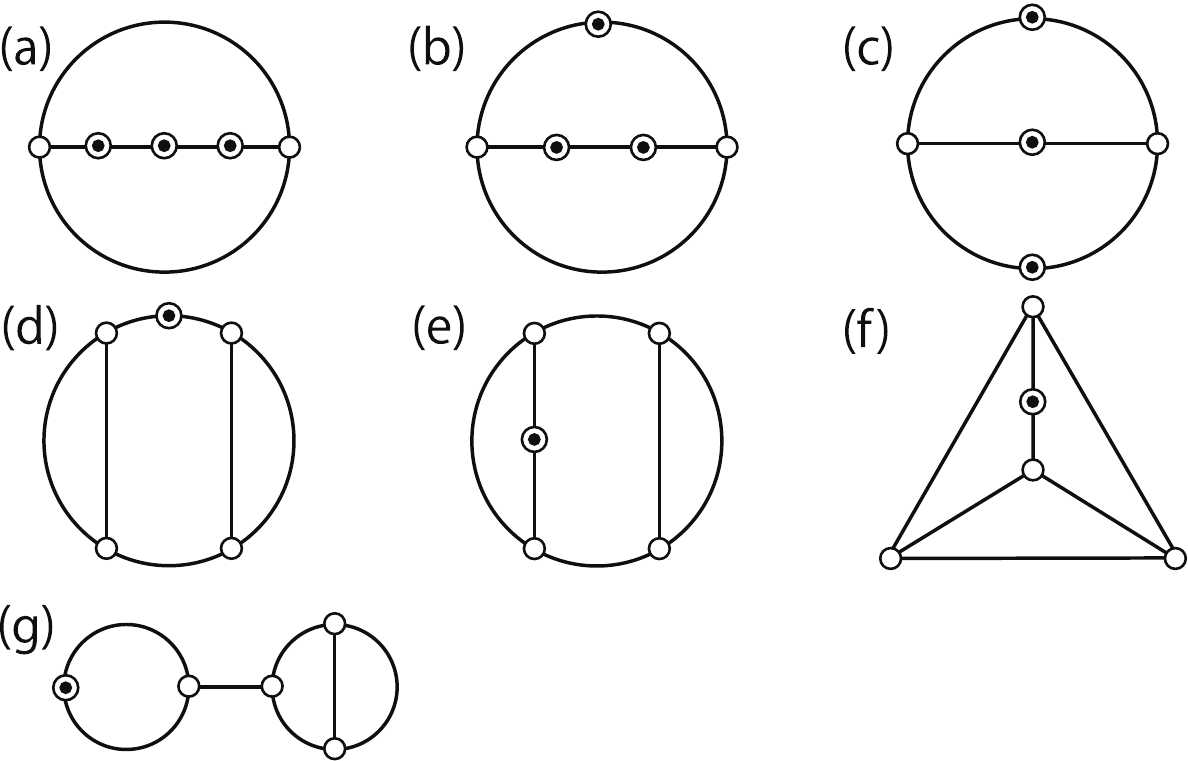}}
\caption{\label{Fig13}
(a),(b),(c) Graphs with three black vertices.
(d),(e),(f),(g) Graphs with one black vertex.
}
\end{figure}

\section{Cases of the graphs as shown in Fig.~\ref{Fig13}(a),(c)}
\label{s:TypeATypeC}

In this section,
we shall show that if $\Gamma$ is a minimal chart of type $(m;5,2)$,
then the graph $\Gamma_m$ contains neither
graphs as shown in Fig.~\ref{Fig13}(a),(c).

\begin{lemma}
\label{LemmaWithTerminal}
{\rm (\cite[Lemma 3.4]{ChartAppVI})}
Let $\Gamma$ be a minimal chart,
and $m$ a label of $\Gamma$.
Let $G$ be a connected component of $\Gamma_m$.
If $G$ contains exactly five white vertices,
and if $G$ has no loop,
then $G$ is one of nine graphs as shown 
in  Fig.~\ref{Fig02} and Fig.~\ref{Fig13}.
\end{lemma}

\begin{figure}[htb]
\centerline{\includegraphics{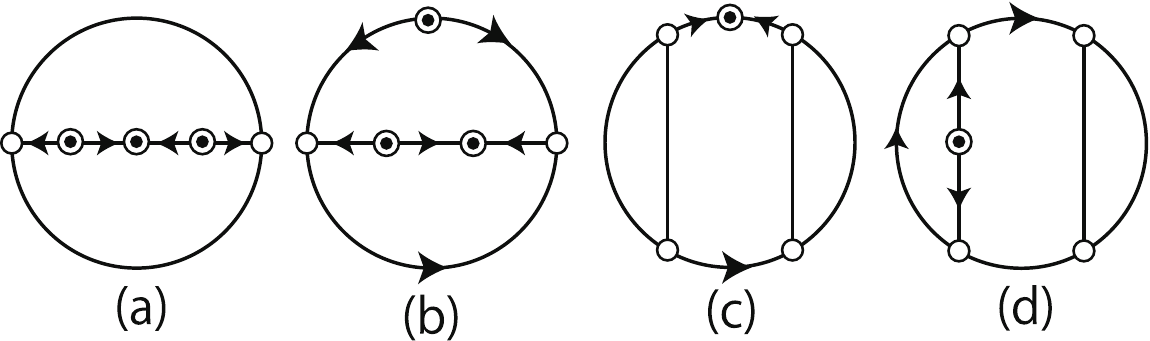}}
\caption{\label{Fig14}
Connected components of $\Gamma_m$ with five white vertices.}
\end{figure}

\begin{lemma}
\label{OriGammaM5}
{\rm (\cite[Lemma 7.2(a),(c)]{ChartAppVIII})}
Let $\Gamma$ be a minimal chart, and $m$ a label of $\Gamma$.
Let $G$ be a connected component of $\Gamma_m$
with $w(G)=5$.
Then we have the following:

\begin{enumerate}
\item[{\rm (a)}]
If $G$ is the graph as shown in Fig.~\ref{Fig13}$($a$)$
$($resp. Fig.~\ref{Fig13}$($b$))$,
then $G$ is one of the RO-family of the graph as shown
in Fig.~\ref{Fig14}$($a$)$
$($resp. Fig.~\ref{Fig14}$($b$))$.
\item[{\rm (b)}]
If $G$ is the graph as shown in Fig.~\ref{Fig13}$($d$)$
$($resp. Fig.~\ref{Fig13}$($e$))$,
then $G$ is one of the RO-family of the graph as shown
in Fig.~\ref{Fig14}$($c$)$
$($resp. Fig.~\ref{Fig14}$($d$))$.
\end{enumerate}
\end{lemma}

\begin{figure}
\centerline{\includegraphics{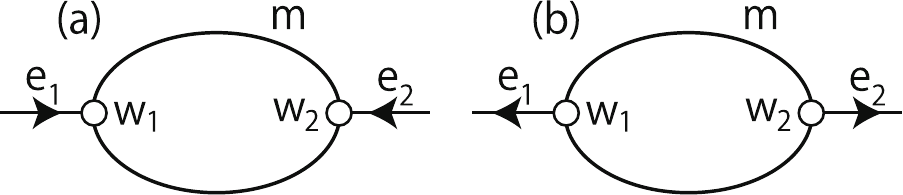}}
\caption{\label{Fig15}
2-angled disks without feelers. }
\end{figure}

\begin{lemma}
\label{2AngledDiskWithoutFeelers}
{\rm (\cite[Lemma 3.6(a)]{ChartAppVIII})}
Let $\Gamma$ be a minimal chart,
and $m$ a label of $\Gamma$.
Let $D$ be a $2$-angled disk of $\Gamma_m$
without feelers, and
$w_1,w_2$ the white vertices in $\partial D$.
Let $e_1,e_2$ be the internal edges
$($possibly terminal edges$)$ of label $m$
at $w_1,w_2$, respectively,
such that $e_1\not\subset D$
and $e_2\not\subset D$.
Suppose that the two edges $e_1,e_2$ are oriented inward
$($resp. outward$)$ at $w_1,w_2$, 
respectively
$($see Fig.~\ref{Fig15}$($a$)$ and $($b$))$.
Then we have 
$w(\Gamma\cap{\rm Int}D)\ge1$.
\end{lemma}

\begin{lemma}
\label{NoTypeA}
Let $\Gamma$ be a minimal chart of type $(m;5,2)$.
Then $\Gamma_m$ does not contain the graph as shown
in Fig.~\ref{Fig13}$($a$)$.
\end{lemma}

\begin{Proof}
Suppose that $\Gamma_m$ contains the graph as shown
in Fig.~\ref{Fig13}$($a$)$, say $G$.
Then $G$ separates the 2-sphere $S^2$ into three disks.
One of the three disks is a 2-angled disk, say $D_1$.
Let $D_2,D_3$ be the other disks.
Then $D_2$ and $D_3$ are special 5-angled disks.
Since one of $D_2,D_3$ has at least two feelers,
by Corollary~\ref{5AngledTwoFeeler}(b)
we have 
\begin{enumerate}
\item[(1)] $w(\Gamma\cap {\rm Int}D_2)\geqq2$ or
$w(\Gamma\cap {\rm Int}D_3)\geqq2$.
\end{enumerate}

By Lemma~\ref{OriGammaM5}(a),
the graph $G$ is one of the RO-family of the graph as shown 
in Fig.~\ref{Fig14}(a).
Without loss of generality,
we can assume that 
the graph $G$ is the graph as shown 
in Fig.~\ref{Fig14}(a). 
Let $w_1,w_2$ be the two white vertices in $\partial D_1$.
Let $e_1,e_2$ be internal edges of label $m$ at $w_1,w_2$,
respectively, with $e_1\not\subset \partial D_1$ and
$e_2\not\subset \partial D_1$.
Then the two edges $e_1,e_2$ are oriented inward at $w_1,w_2$,
respectively.
Thus by Lemma~\ref{2AngledDiskWithoutFeelers},
we have $w(\Gamma\cap{\rm Int}D_1)\geqq1$.
Hence by (1), we have

$\begin{array}{rcl}
7=w(\Gamma) & = & w(G)+w(\Gamma\cap{\rm Int}D_1)+
w(\Gamma\cap{\rm Int}D_2)+w(\Gamma\cap{\rm Int}D_3)\vspace{2mm}\\
& \geqq & 5+1+2=8.\vspace{2mm}
\end{array}
$\\
This is a contradiction.
Therefore $\Gamma_m$ does not contain the graph as shown
in Fig.~\ref{Fig13}$($a$)$.
\end{Proof}

\begin{lemma}
\label{NoTypeC}
Let $\Gamma$ be a minimal chart of type $(m;5,2)$.
Then $\Gamma_m$ does not contain the graph as shown
in Fig.~\ref{Fig13}$($c$)$.
\end{lemma}

\begin{Proof}
Suppose that $\Gamma_m$ contains the graph as shown
in Fig.~\ref{Fig13}(c), say $G$.
Then $G$ separates the 2-sphere $S^2$ into three disks.
Let $D_1,D_2,D_3$ be the three disks.
Then $D_1,D_2,D_3$ are special 4-angled disks.
Without loss of generality,
we can assume that $D_1$ has at least one feeler.
Then $D_1$ has one feeler or two feelers.

If $D_1$ has two feelers,
then one of $D_2,D_3$ has one feeler.
Thus by Colorally~\ref{4AngledTwoFeeler},
we have $w(\Gamma\cap{\rm Int}D_1)\geqq2$ and
($w(\Gamma\cap{\rm Int}D_2)\geqq 1$ or 
$w(\Gamma\cap{\rm Int}D_3)\geqq 1$).
By the similar way of the proof of Lemma~\ref{NoTypeA},
we have a contradiction.
Thus $D_1$ has exactly one feeler.

Since $D_1$ has exactly one feeler,
one of $D_2,D_3$ has at least one feeler.
If one of $D_2,D_3$ has exactly two feelers,
then we have the same contradiction as above.
Hence both of $D_2,D_3$ have exactly one feeler.
Thus by Colorally~\ref{4AngledTwoFeeler}(a),
we have $w(\Gamma\cap{\rm Int}D_i)\geqq1$ for $i=1,2,3$.
Hence\vspace{2mm}

$\begin{array}{rcl}
7=w(\Gamma) & = & w(G)+w(\Gamma\cap{\rm Int}D_1)+
w(\Gamma\cap{\rm Int}D_2)+w(\Gamma\cap{\rm Int}D_3)\vspace{2mm}\\
& \geqq & 5+1+1+1=8.\vspace{2mm}
\end{array}
$\\
This is a contradiction.
Therefore $\Gamma_m$ does not contain the graph as shown
in Fig.~\ref{Fig13}(c).
\end{Proof}



\section{Case of the graph as shown in Fig.~\ref{Fig13}(b)}
\label{s:TypeB}

In this section,
we shall show that if $\Gamma$ is a minimal chart of type $(m;5,2)$,
then the graph $\Gamma_m$ does not contain 
the graph as shown in Fig.~\ref{Fig13}(b).

\begin{figure}[htb]
\centerline{\includegraphics{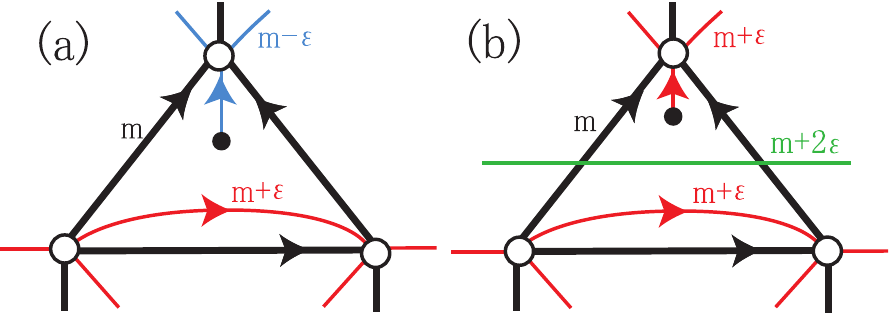}}
\caption{\label{Fig16}
The 3-angled disks  have no feelers,
$m$ is a label, $\varepsilon\in\{+1,-1\}$.}
\end{figure}

\begin{lemma}
{\rm (\cite[Lemma 4.2(a)]{ChartAppIX})}
\label{Theorem3AngledDisk}
Let $\Gamma$ be a minimal chart, and $m$ a label of $\Gamma$.
Let $D$ be a special $3$-angled disk of $\Gamma_m$
with at most two feelers.
If $w(\Gamma\cap {\rm Int}D)=0$,
then a regular neighborhood of $D$ contains one of the RO-families of the two pseudo charts as shown in 
Fig.~\ref{Fig16}.

\end{lemma}

\begin{lemma}
\label{NoTypeB}
Let $\Gamma$ be a minimal chart of type $(m;5,2)$.
Then $\Gamma_m$ does not contain the graph as shown
in Fig.~\ref{Fig13}$($b$)$.
\end{lemma}

\begin{figure}[htb]
\centerline{\includegraphics{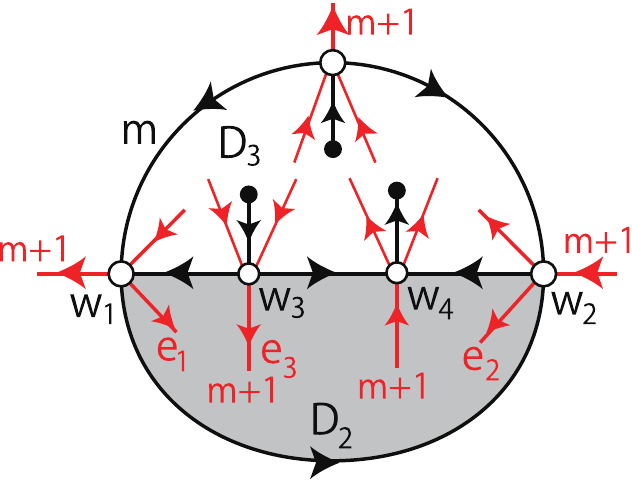}}
\caption{\label{Fig17}
The graph as shown in Fig.~\ref{Fig13}(b).
The gray region is the 4-angled disk $D_2$.}
\end{figure}

\begin{Proof}
Suppose that $\Gamma_m$ contains the graph as shown
in Fig.~\ref{Fig13}(b), say $G$.
Then $G$ separates the 2-sphere $S^2$ into three disks.
One of the three disks is a 3-angled disk, say $D_1$.
One of the three disks is a 4-angled disk, say $D_2$.
The last one is a 5-angled disk, say $D_3$.
For the disk $D_3$,
there are four cases:
(i) $D_3$ has no feeler,
(ii) $D_3$ has exactly one feeler,
(iii) $D_3$ has exactly two feelers,
(iv) $D_3$ has exactly three feelers.

{\bf Case (i).}
Since $D_3$ has no feeler,
the 3-angled disk $D_1$ has exactly one feeler and
the 4-angled disk $D_2$ has exactly two feelers.
Thus by Corollary~\ref{4AngledTwoFeeler}(b)
and Lemma~\ref{Theorem3AngledDisk},
we have  $w(\Gamma\cap {\rm Int}D_1)\geqq1$ and
$w(\Gamma\cap {\rm Int}D_2)\geqq2$.
Hence, we have

$\begin{array}{rcl}
7=w(\Gamma) & \geqq & w(G)+w(\Gamma\cap{\rm Int}D_1)+
w(\Gamma\cap{\rm Int}D_2)\vspace{2mm}\\
& \geqq & 5+1+2=8.\vspace{2mm}
\end{array}
$\\
This is a contradiction.
Thus Case (i) does not occur.

{\bf Case (ii).}
Since the 5-angled disk $D_3$ has exactly one feeler,
by Corollary~\ref{5AngledTwoFeeler}(a)
we have 
$w(\Gamma\cap {\rm Int}D_3)\geqq 1$.
Moreover,
the 4-angled disk $D_2$ has one feeler or two feelers.

If $D_2$ has exactly two feelers,
then by Corollary~\ref{4AngledTwoFeeler}(b)
we have $w(\Gamma\cap {\rm Int}D_2)\geqq2$.
Thus we have the same contradiction of Case (i).

If $D_2$ has exactly one feeler,
then the 3-angled disk $D_1$ has exactly one feeler.
Thus by Corollary~\ref{4AngledTwoFeeler}(a) and Lemma~\ref{Theorem3AngledDisk},
we have  $w(\Gamma\cap {\rm Int}D_2)\geqq1$ and
$w(\Gamma\cap {\rm Int}D_3)\geqq1$.
By the similar way of the proof of Lemma~\ref{NoTypeC},
we have a contradiction.

Therefore both cases do not occur. Thus Case (ii) does not occur.

{\bf Case (iii).}
Since the 5-angled disk $D_3$ has exactly two feelers,
one of the disks $D_1,D_2$ has exactly one feeler.
Thus by Corollary~\ref{5AngledTwoFeeler}(b),
Corollary~\ref{4AngledTwoFeeler}(a)
and Lemma~\ref{Theorem3AngledDisk},
we have $w(\Gamma\cap {\rm Int}D_3)\geqq2$ and
($w(\Gamma\cap {\rm Int}D_1)\geqq1$ or 
$w(\Gamma\cap {\rm Int}D_2)\geqq1$).
Hence we have the same contradiction of Case (i).
Thus Case (iii) does not occur.

{\bf Case (iv).}
Since the 5-angled disk $D_3$ has exactly three feelers,
by Corollary~\ref{5AngledTwoFeeler}(b)
we have $w(\Gamma\cap {\rm Int}D_3)\geqq2$.
Since $w(\Gamma)=7$,
we have $w(\Gamma\cap {\rm Int}D_3)=2$ and $w(\Gamma\cap {\rm Int}D_2)=0$.

By Lemma~\ref{OriGammaM5}(a),
the graph $G$ is one of the RO-family of the graph as shown 
in Fig.~\ref{Fig14}(b).
Without loss of generality,
we can assume that 
the graph $G$ is the graph as shown 
in Fig.~\ref{Fig14}(b). 
Thus the chart $\Gamma$ contains the pseudo chart
as shown in Fig.~\ref{Fig17}.
We use the notations as shown in Fig.~\ref{Fig17},
where $e_1,e_2,e_3$ are internal edges 
(possibly terminal edges) of label $m+1$
oriented outward at $w_1,w_2,w_3$ in $D_2$,
respectively.
Hence the condition $w(\Gamma\cap{\rm Int}D_2)=0$
implies
one of $e_1$ or $e_2$ is a terminal edge.
However neither $e_1$ nor $e_2$ is
middle at $w_1$ or $w_2$
(see Fig.~\ref{Fig17}).
This contradicts Assumption~\ref{AssumeTerminal}.
Thus Case (iv) does not occur.

Therefore all four cases do not occur.
Hence $\Gamma_m$ does not contain the graph as shown
in Fig.~\ref{Fig13}(b).
\end{Proof}

\section{Case of the graph as shown in Fig.~\ref{Fig13}(d)}
\label{s:TypeD}

In this section,
we shall show that if $\Gamma$ is a minimal chart of type $(m;5,2)$,
then the graph $\Gamma_m$ does not contain
the graph as shown in Fig.~\ref{Fig13}(d).

\begin{lemma}
\label{NoTypeD}
Let $\Gamma$ be a minimal chart of type $(m;5,2)$.
Then $\Gamma_m$ does not contain the graph as shown
in Fig.~\ref{Fig13}$($d$)$.
\end{lemma}

\begin{figure}[htb]
\centerline{\includegraphics{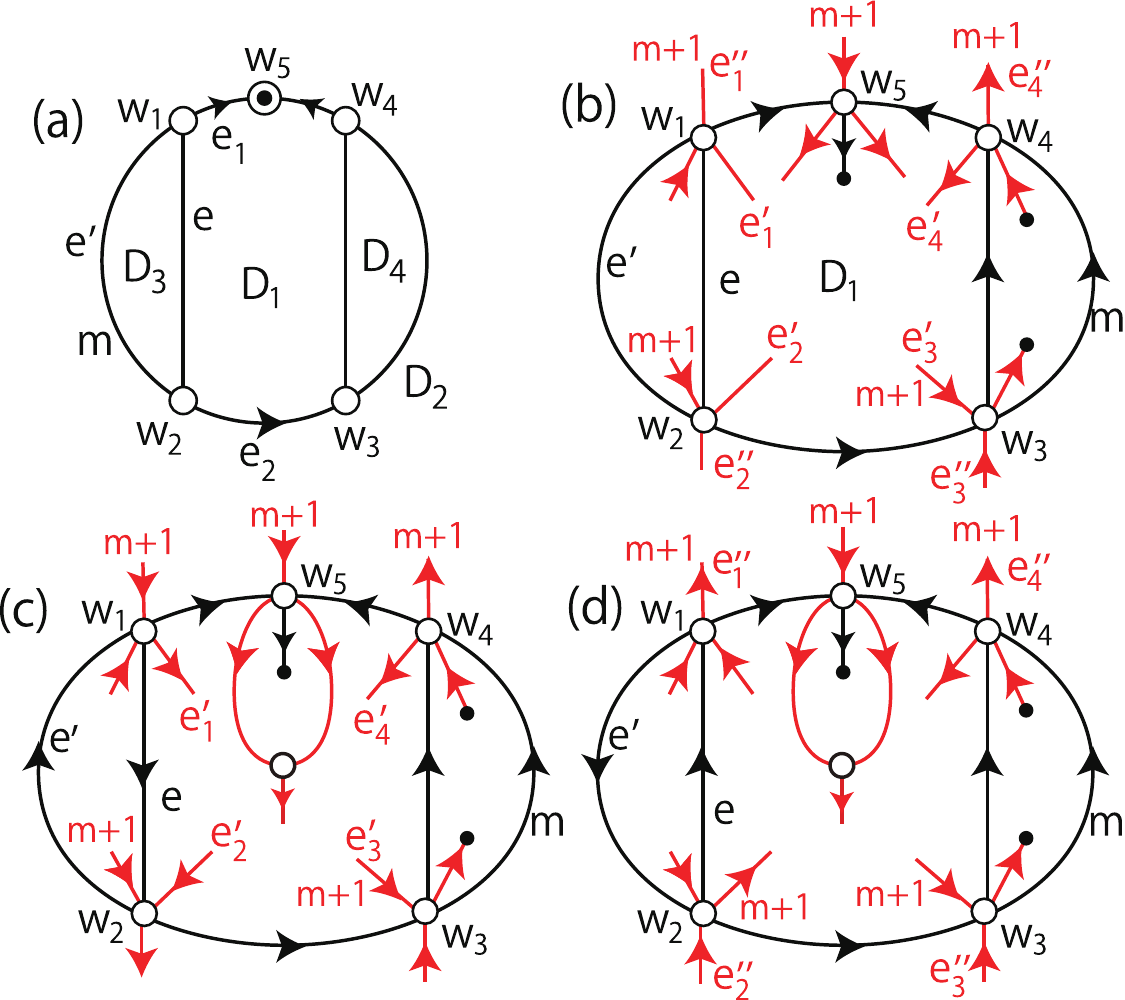}}
\caption{\label{Fig18}
The graphs as shown in Fig.~\ref{Fig13}(d).
(c) The edge $e$ is oriented from $w_1$ to $w_2$.
(d) The edge $e$ is oriented from $w_2$ to $w_1$.}
\end{figure}

\begin{Proof}
Suppose that $\Gamma_m$ contains the graph as shown
in Fig.~\ref{Fig13}(d), say $G$.
By Lemma~\ref{OriGammaM5}(b),
the graph $G$ is one of the RO-family of the graph as shown 
in Fig.~\ref{Fig14}(c).
Without loss of generality,
we can assume that 
the graph $G$ is the graph as shown 
in Fig.~\ref{Fig14}(c). 
We use the notations as shown in Fig.~\ref{Fig18}(a),
where  $w_1,w_2,\cdots,w_5$ are five white vertices, and
\begin{enumerate}
\item[(1)] $e_1,e_2$ are internal edges of label $m$
oriented outward at $w_1,w_2$, respectively.
\end{enumerate}

Let $D_1,D_2$ be special 5-angled disks of $\Gamma_m$
with ${\rm Int}D_1\cap{\rm Int}D_2=\emptyset$ such that
the disk $D_2$ contains the point at infinity, $\infty$.
If necessary, we move the point $\infty$ by Assumption~\ref{Infinity}, and if necessary, we reflect the chart $\Gamma$,
we can assume that the disk $D_1$ has one feeler.
Thus by Corollary~\ref{5AngledTwoFeeler}(a),
we have 
\begin{enumerate}
\item[(2)] $w(\Gamma\cap{\rm Int}D_1)\geqq 1$.
\end{enumerate}

Let $D_3,D_4$ be special 2-angled disks of $\Gamma_m$
with $\partial D_3\ni w_1$ and $\partial D_4\ni w_4$.
Then by (1) and Lemma~\ref{2AngledDiskWithoutFeelers},
we have $w(\Gamma\cap{\rm Int}D_3)\geqq 1$.
Hence by (2), the condition $w(\Gamma)=7$ implies that
\begin{enumerate}
\item[(3)] $w(\Gamma\cap{\rm Int}D_1)= 1$, 
$w(\Gamma\cap{\rm Int}D_2)= 0$,
$w(\Gamma\cap{\rm Int}D_4)= 0$.
\end{enumerate}
Thus by Lemma~\ref{Theorem2AngledDisk},
a regular neighborhood of $D_4$ contains
the pseudo chart as shown in Fig.~\ref{Fig10}(b).
Therefore, the chart $\Gamma$ contains
the pseudo chart as shown in Fig.~\ref{Fig18}(b).
We use the notations as shown in Fig.~\ref{Fig18}(b),
where $e_i',e_i''$ ($i=1,2,3,4$) are internal edges
(possibly terminal edges) of label $m+1$ at $w_i$
with $e_i'\subset D_1$ and $e_i''\subset D_2$,
\begin{enumerate}
\item[(4)] $e_3',e_3''$ are oriented inward at $w_3$,
but not middle at $w_3$, and
\end{enumerate}
neither $e_4'$ nor $e_4''$ is middle at $w_4$.
Thus by Assumption~\ref{AssumeTerminal},
\begin{enumerate}
\item[(5)]
none of $e_3',e_3'',e_4',e_4''$ are terminal edges.
\end{enumerate}

Since $w(\Gamma\cap{\rm Int}D_1)=1$ by (3) and
since $e_4'$ is not a terminal edge by (5),
by Lemma~\ref{5AngledOneFeelerOneWhite}
the disk $D_1$ contains the pseudo chart
as shown in Fig.~\ref{Fig11}(e)
(see Fig.~\ref{Fig18}(c),(d)).
Let $e,e'$ be internal edges of label $m$ in $\partial D_3$
with $e\subset D_1$ and $e'\subset D_2$.
Then there are two cases:
(i) $e$ is oriented from $w_1$ to $w_2$ 
(see Fig.~\ref{Fig18}(c)),
(ii) $e$ is oriented from $w_2$ to $w_1$ 
(see Fig.~\ref{Fig18}(d)).

{\bf Case (i).}
By looking around the white vertex $w_1$,
we have that $e'$ is oriented from $w_2$ to $w_1$, and
\begin{enumerate}
\item[(6)]
$e_2'$ is oriented inward at $w_2$, but not middle at $w_2$.
\end{enumerate}
Thus by (5) and Assumption~\ref{AssumeTerminal},
none of $e_2',e_3',e_4'$ are terminal edges.
Moreover, by (4) and (6),
the two edges $e_2'$ and $e_3'$ are oriented inward at $w_2,w_3$,
respectively.
Hence, for the edge $e_3'$, we must have $e_3'=e_4'$.
However, there exists a lens.
This contradicts Lemma~\ref{NoLens}.
Thus Case (i) does not occur.

{\bf Case (ii).}
Looking around the white vertex $w_2$,
we have that $e'$ is oriented from $w_1$ to $w_2$, and
\begin{enumerate}
\item[(7)]
$e_2''$ is oriented inward at $w_2$, but not middle at $w_2$.
\end{enumerate}
Now $w(\Gamma\cap{\rm Int}D_2)=0$ by (3).
By the similar way of the proof of Case (i) in this lemma,
for the edge $e_3''$, we must have $e_3''=e_4''$.
However, there exists a lens.
This contradicts Lemma~\ref{NoLens}.
Thus Case (ii) does not occur.

Therefore both two cases do not occur.
Hence $\Gamma_m$ does not contain the graph as shown
in Fig.~\ref{Fig13}(d).
\end{Proof}

\section{IO-Calculation}

\label{s:IOC}

In this section,
we review IO-Calculation.

Let $\Gamma$ be a chart,
 and $v$ a vertex. 
Let $\alpha$ be a short arc of $\Gamma$ in a small neighborhood of $v$ such that $v$ is an endpoint of $\alpha$. 
If the arc $\alpha$ is oriented to $v$, then $\alpha$ is called {\it an inward arc}, 
and otherwise $\alpha$ is called {\it an outward arc}.

\begin{figure}
\centerline{\includegraphics{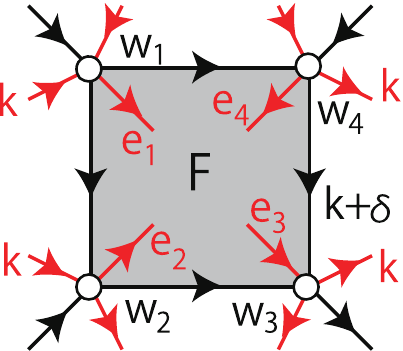}}
\caption{\label{Fig19} The gray region is the 4-angled disk $F$, $k$ is a label, $\delta\in\{+1,-1\}$.}
\end{figure}

Let $\Gamma$ be an $n$-chart. 
Let $F$ be a closed domain with $\partial F\subset \Gamma_{k-1}\cup\Gamma_{k}\cup \Gamma_{k+1}$ for some label $k$ of $\Gamma$, where $\Gamma_0=\emptyset$ and $\Gamma_{n}=\emptyset$. 
By Condition (iii) for charts,
in a small neighborhood of each white vertex, there are three inward arcs and three outward arcs.
Also in a small neighborhood of each black vertex, there exists only one inward arc or one outward arc.
We often use the following fact, 
when we fix (inward or outward) arcs 
near white vertices and black vertices: 
\begin{enumerate}
\item[$(*)$]
{\it The number of inward arcs contained in $F\cap \Gamma_k$ is equal to the number of outward arcs in $F\cap \Gamma_k$.
}
\end{enumerate}
When we use this fact, 
we say that we use {\it IO-Calculation with respect to $\Gamma_k$ in $F$}.
For example, in a minimal chart $\Gamma$, 
consider the pseudo chart as shown in Fig.~\ref{Fig19} 
where
\begin{enumerate}
\item[(1)] $F$ is a $4$-angled disk of $\Gamma_{k+\delta}$ 
without feelers for some $\delta\in\{+1,-1\}$,
\item[(2)] $e_1,e_2,e_4$ are internal edges (possibly terminal edges) of label $k$ oriented outward at $w_1,w_2,w_4$, respectively, \item[(3)]  $e_3$ is an internal edge (possibly a terminal edge) 
of label $k$ oriented inward at $w_3$, 
\item[(4)] neither $e_2$ nor $e_4$ is middle at $w_2$ or $w_4$.
\end{enumerate}
Then we can show that $w(\Gamma\cap{\rm Int}F)\ge1$.
Suppose $w(\Gamma\cap{\rm Int}F)=0$.
By (4) and Assumption~\ref{AssumeTerminal},
\begin{enumerate}
\item[(5)] neither $e_2$ nor $e_4$ is a terminal edge. 
\end{enumerate}

If both two edges $e_1,e_3$ are a terminal edge,
then 
by (2) and (3) 
the number of inward arcs in $F\cap \Gamma_k$ is two,  
but the number of outward arcs in $F\cap \Gamma_k$ is four. 
This contradicts the fact $(*)$. 
If $e_1$ is a terminal edge,
but $e_3$ is not a terminal edge,
then
by (2) and (3) 
the number of inward arcs in $F\cap \Gamma_k$ is two,  
but the number of outward arcs in $F\cap \Gamma_k$ is three. 
This contradicts the fact~$(*)$. 
Similarly for the other cases
 we have the same contradiction.
Thus $w(\Gamma\cap{\rm Int}F)\ge1$.
Instead of the above argument, 
\begin{enumerate}
\item[]
{\it we have $w(\Gamma\cap{\rm Int}F)\ge1$ 
by IO-Calculation with respect to $\Gamma_{k}$ in $F$.}
\end{enumerate}

\section{Case of the graph as shown in Fig.~\ref{Fig13}(e)}
\label{s:TypeE}

In this section,
we shall show that if $\Gamma$ is a minimal chart of type $(m;5,2)$,
then the graph $\Gamma_m$ does not contain
the graph as shown in Fig.~\ref{Fig13}(e).

\begin{figure}[htb]
\begin{center}
\centerline{\includegraphics{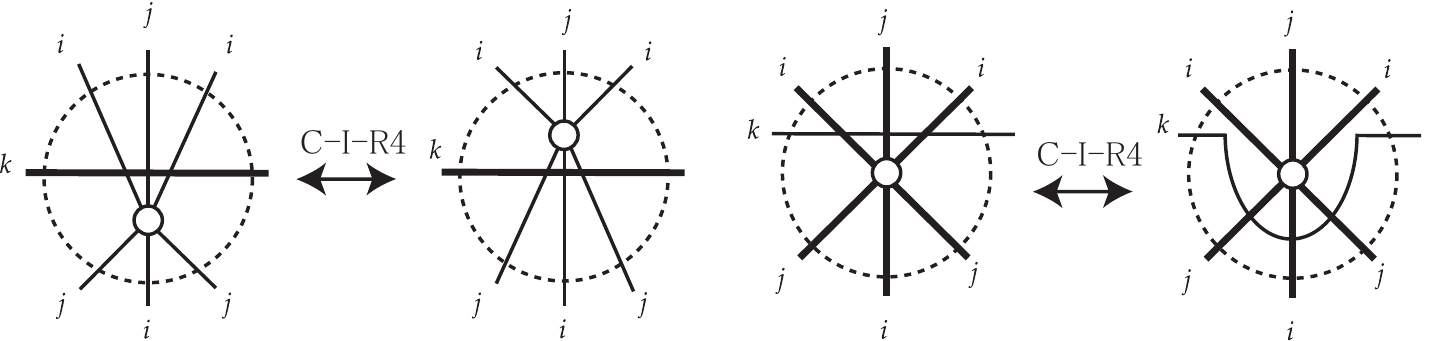}}
\caption{\label{Fig20}
C-moves keeping thicken figures fixed.}
\end{center}
\end{figure}

Let $\Gamma $ and $\Gamma^\prime $ be C-move equivalent charts. 
Suppose that a pseudo chart $X$ of $\Gamma$ is also a pseudo chart of $\Gamma^\prime$. 
Then we say that 
$\Gamma$ is modified to $\Gamma^\prime$ by {\it C-moves keeping $X$ fixed}.
In Fig.~\ref{Fig20},
we give examples of C-moves keeping pseudo charts  fixed.

Let $\Gamma$ be a chart, and $X$ a subset of $\Gamma$. Let 
$$c(X)=\text{the number of crossings in $X$}.$$

Let $D$ be a $k$-angled disk of $\Gamma_m$ for a minimal chart $\Gamma$.
The pair of integers 
$(w(\Gamma\cap {\rm Int}D),c(\partial D))$
is called the {\it local complexity 
with respect to $D$},
denoted by $\ell c(D;\Gamma)$.
Let
${\Bbb S}$ be the set of all minimal charts each of which can be moved from $\Gamma$ by C-moves in a regular neighborhood of $D$ keeping $\partial D$ fixed.
The chart $\Gamma$ is said to be 
{\it locally minimal
with respect to $D$}
if its local complexity
with respect to $D$
is minimal
among the charts in ${\Bbb S}$ with respect to 
the lexicographic order.

\begin{figure}[htb]
\begin{center}
\centerline{\includegraphics{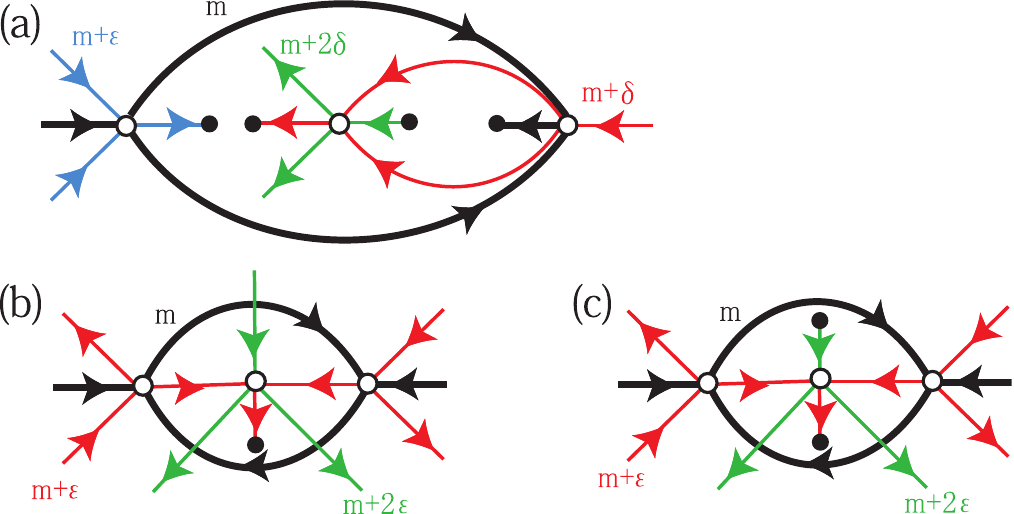}}
\caption{\label{Fig21}
The 2-angled disk (a) has one feeler, the others do not have any feelers.} 
\end{center}
\end{figure}

\begin{lemma}
\label{Theorem2AngledDisk2}
{\rm (\cite[Theorem 1.1]{ChartAppIII})}
Let $\Gamma$ be a minimal chart.
Let $D$ be a $2$-angled disk of $\Gamma_m$ with at most one feeler
such that $\Gamma$ is locally minimal
with respect to $D$.
If $w(\Gamma\cap{\rm Int}D)\leqq1$,
then a regular neighborhood of $D$ contains an element in the RO-families of the five pseudo charts as shown in 
Fig.~\ref{Fig10} and Fig.~\ref{Fig21}.
\end{lemma}

\begin{lemma}
\label{NoTypeE}
Let $\Gamma$ be a minimal chart of type $(m;5,2)$.
Then $\Gamma_m$ does not contain the graph as shown
in Fig.~\ref{Fig13}$($e$)$.
\end{lemma}

\begin{figure}[htb]
\centerline{\includegraphics{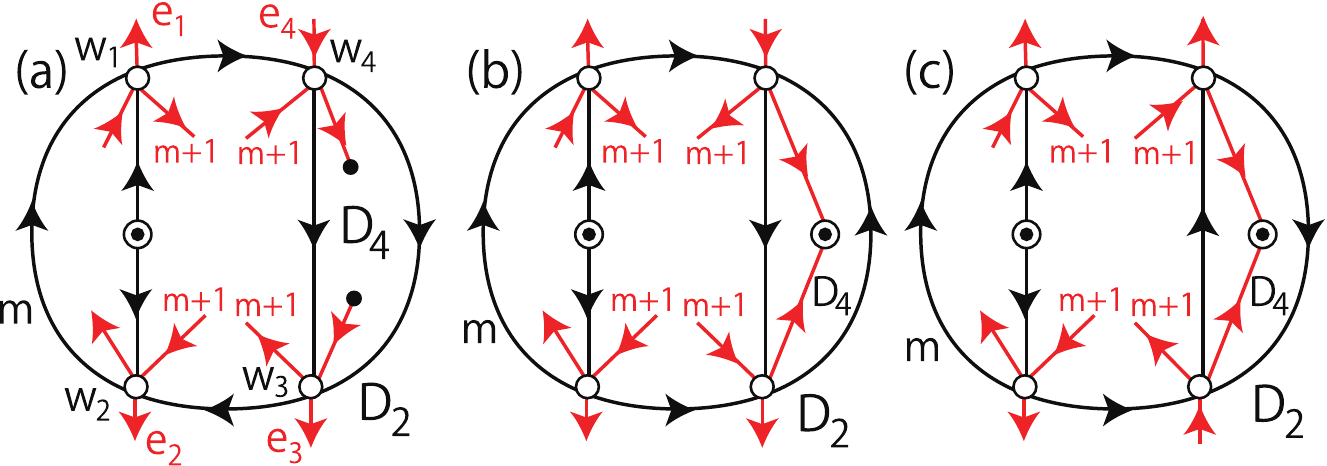}}
\caption{\label{Fig22}
The graphs as shown in Fig.~\ref{Fig13}(e).
(a) The two internal edges in $\partial D_4$ are oriented from $w_4$ to $w_3$.
(b) $\partial D_4$ is oriented anticlockwise.
(c) $\partial D_4$ is oriented clockwise.}
\end{figure}

\begin{Proof}
Suppose that $\Gamma_m$ contains the graph as shown
in Fig.~\ref{Fig13}(e), say $G$.
Then $G$ separates the 2-sphere $S^2$ into
four disks.
One of the four disks is a 5-angled disk, say $D_1$.
One of the four disks is a 4-angled disk, say $D_2$.
One of the four disks is a 3-angled disk, say $D_3$.
Let $D_4$ be the last disk.

Since one of $D_1,D_3$ has one feeler,
by Corollary~\ref{5AngledTwoFeeler}(a) and 
Lemma~\ref{Theorem3AngledDisk}
we have
\begin{enumerate}
\item[(1)] $w(\Gamma\cap{\rm Int}D_1)\geqq 1$ or
$w(\Gamma\cap{\rm Int}D_3)\geqq 1$.
\end{enumerate}
Thus the condition $w(\Gamma)=7$ implies that
$w(\Gamma\cap{\rm Int}D_4)\leqq 1$.
We can assume that $\Gamma$ is locally minimal with
respect to $D_4$.
Hence by Lemma~\ref{Theorem2AngledDisk2}
a regular neighborhood of $D_4$
contains one of the RO-families of
the three pseudo charts as shown in
Fig.~\ref{Fig10}(b) and Fig.~\ref{Fig21}(b),(c).
Moreover, by Lemma~\ref{OriGammaM5}(b),
the graph $G$ is the graph as shown in 
Fig.~\ref{Fig14}(d).
Thus the chart $\Gamma$ 
contains one of the RO-families
of the three pseudo charts as shown in 
Fig.~\ref{Fig22},
where 
the pseudo charts as shown in Fig.~\ref{Fig22}(b),(c)
are contained in one of the pseudo charts as shown in
Fig.~\ref{Fig21}(b),(c).
Without loss of generality,
we can assume that 
the chart $\Gamma$ contains one of 
the three pseudo charts as shown in 
Fig.~\ref{Fig22}.

Suppose that 
the chart $\Gamma$ contains the pseudo chart as shown in 
Fig.~\ref{Fig22}(b).
Then we have $w(\Gamma\cap{\rm Int}D_4)\geqq1$.
Moreover, we have $w(\Gamma\cap{\rm Int}D_2)\geqq1$
by considering as $F = D_2$, $k= m+1$ and 
$\delta=-1$ in the example
of IO-Calculation in Section~\ref{s:IOC}.
Hence by (1)\vspace{2mm}\\
$\begin{array}{rl}
7=& w(\Gamma)\vspace{2mm}\\
= & w(G)+w(\Gamma\cap{\rm Int}D_1)+
w(\Gamma\cap{\rm Int}D_2)+w(\Gamma\cap{\rm Int}D_3)+
w(\Gamma\cap{\rm Int}D_4)\vspace{2mm}\\
\geqq & 5+1+1+1=8.\vspace{2mm}
\end{array}
$\\
This is a contradiction.
Thus $\Gamma$ does not contain the pseudo chart as shown in 
Fig.~\ref{Fig22}(b).

Similarly, we can show that
 $\Gamma$ does not contain the pseudo chart as shown in 
Fig.~\ref{Fig22}(c).

Now, suppose that 
the chart $\Gamma$ contain the pseudo chart as shown in 
Fig.~\ref{Fig22}(a).
We use the notations as shown in 
Fig.~\ref{Fig22}(a), where
$e_i$ $(i=1,2,3,4)$ is an internal edge
(possibly a terminal edge) of label $m+1$ at $w_i$
in the 4-angled disk $D_2$, and
\begin{enumerate}
\item[(2)] the three edges $e_1,e_2,e_3$ are oriented
outward at $w_1,w_2,w_3$, respectively,
\end{enumerate}
none of  $e_1,e_2,e_3,e_4$ are middle at $w_1,w_2,w_3$ or $w_4$.
Thus by Assumption~\ref{AssumeTerminal},
\begin{enumerate}
\item[(3)] none of  $e_1,e_2,e_3,e_4$ are terminal edges.
\end{enumerate}
Hence by IO-Calculation with respect to 
$\Gamma_{m+1}$ in $D_2$,
we have $w(\Gamma\cap{\rm Int}D_2)\geqq1$.
Thus by (1),
the condition $w(\Gamma)=7$ implies that
\begin{enumerate}
\item[(4)] $w(\Gamma\cap{\rm Int}D_2)=1$.
\end{enumerate}

Let $w_5$ be the white vertex in ${\rm Int}D_2$.
Then for the edge $e_4$,
there are four cases:
(i) $e_4=e_1$,
(ii) $e_4=e_2$,
(iii) $e_4=e_3$,
(iv) $e_4\ni w_5$.

{\bf Case (i) and Case (iii).}
There exists a lens.
This contradicts Lemma~\ref{NoLens}.
Hence neither Case~(i) nor Case~(iii) occurs.

{\bf Case (ii).}
By (2) and (3),
both of $e_1,e_3$ contain the white vertex $w_5$.
Thus one of the edges $e_1,e_3$ of label $m+1$
intersects the edge $e_4$ of label $m+1$.
This contradicts the definition of the chart.
Hence Case (ii) does not occur.

{\bf Case (iv).}
Since $e_4$ contains the white vertex $w_5$,
one of the three edges $e_1,e_2,e_3$
does not contain white vertex $w_5$.
Thus by (4),
one of the three edges $e_1,e_2,e_3$ is a terminal edge.
This contradicts (3).
Hence Case (iv) does not occur.

Therefore all the four cases do not occur.
Thus $\Gamma$ does not contain the pseudo chart as shown in 
Fig.~\ref{Fig22}(a).
Hence $\Gamma_m$ does not contain the graph as shown
in Fig.~\ref{Fig13}(e).
\end{Proof}

\section{Case of the graph as shown in Fig.~\ref{Fig13}(f)}
\label{s:TypeF}

In this section,
we shall show that if $\Gamma$ is a minimal chart of type $(m;5,2)$,
then the graph $\Gamma_m$ does not contain
the graph as shown in Fig.~\ref{Fig13}(f).

\begin{figure}[htb]
\centerline{\includegraphics{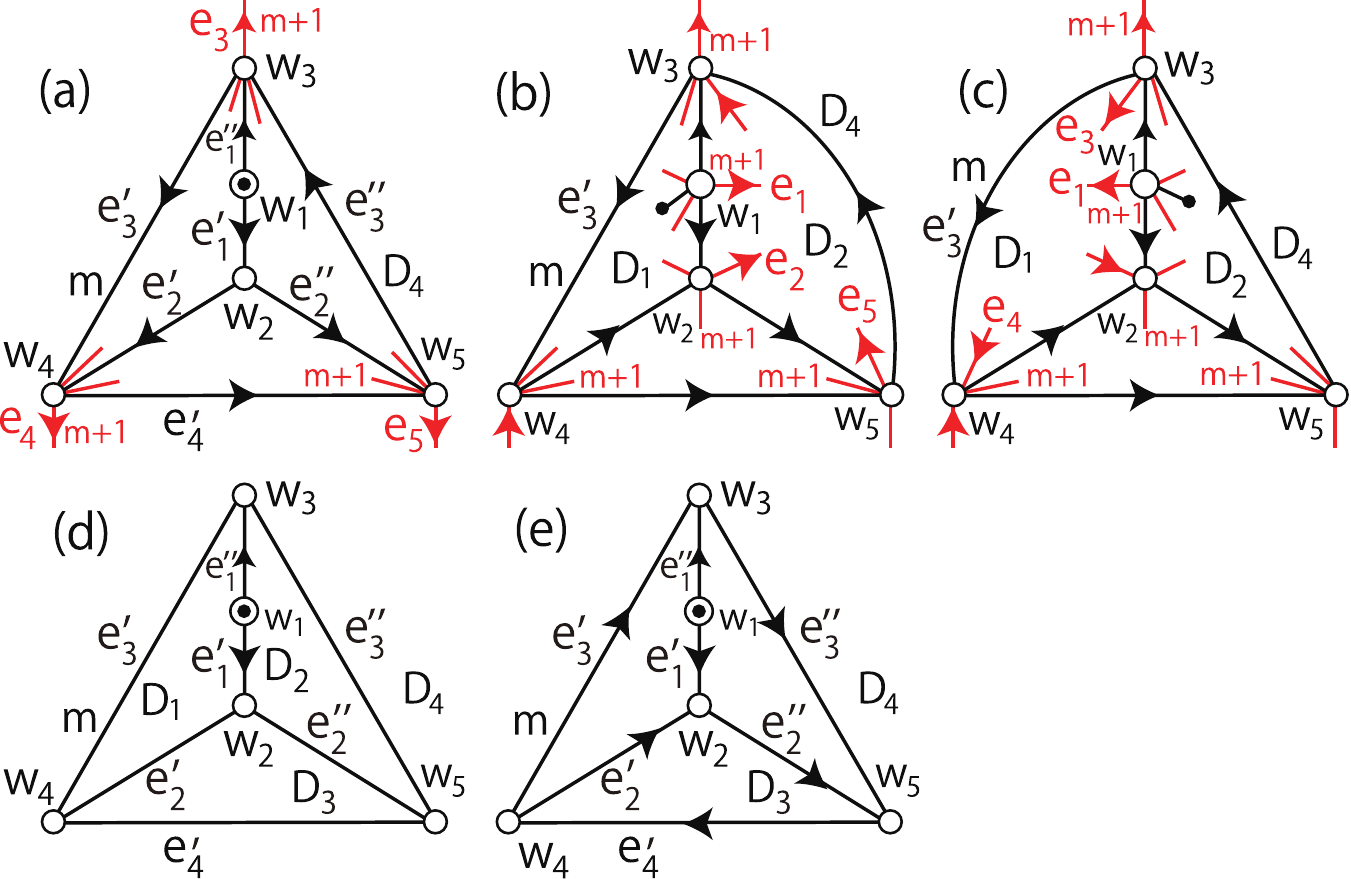}}
\caption{\label{Fig23}
The graphs as shown in Fig.~\ref{Fig13}(f).}
\end{figure}

\begin{lemma}
\label{OriNoTypeF}
Let $\Gamma$ be a minimal chart of type $(m;5,2)$.
If $\Gamma_m$ contains the graph as shown
in Fig.~\ref{Fig13}$($f$)$,
then $\Gamma$ contains one of RO-families
of the three pseudo charts as shown in 
Fig~\ref{Fig23}$($a$)$,$($b$)$,$($c$)$.
\end{lemma}

\begin{Proof}
Let $G$ be the graph in $\Gamma_m$ 
as shown in Fig.~\ref{Fig13}(f).
We use the notations as shown in 
Fig.~\ref{Fig23}(d),
where $w_1$ is the BW-vertex,
and $e_1',e_1'',e_2',e_2'',e_3',e_3'',e_4'$ are 
seven internal edges of
label $m$ with
$e_1'\cap e_1''\ni w_1$,
$e_1'\cap e_2'\cap e_2''\ni w_2$,
$e_1''\cap e_3'\cap e_3''\ni w_3$, and
$\partial e_4'=\{w_4,w_5\}$.

Since the graph $G$ separates the 2-sphere $S^2$ 
into four disks.
Two of the four disks are 4-angled disks,
say $D_1,D_2$.
Two of the four disks are 3-angled disks,
say $D_3,D_4$.
Without loss of generality we can assume that
$\partial D_1\ni w_4$, 
$\partial D_2\ni w_5$
and
$D_4$ contains the point at infinity, $\infty$
(see Fig.~\ref{Fig23}(d)).

Without loss of generality
we can assume that
the terminal edge of label $m$ at $w_1$
is oriented inward at $w_1$. 
Then by Assumption~\ref{AssumeTerminal},
\begin{enumerate}
\item[(1)] both of $e_1',e_1''$ are oriented outward at $w_1$ 
(see Fig.~\ref{Fig23}(d)). 
\end{enumerate}
There are two cases:
(i) one of $e_1',e_1''$ is middle at $w_2$ or $w_3$,
(ii) neither $e_1'$ nor $e_1''$ is middle at $w_2$ or $w_3$.

{\bf Case (i).}
If necessary we move the point  $\infty$ in $D_3$,
we can assume that $e_1'$ is middle at $w_2$.
By Condition~(iii) of the definition of a chart,
both of $e_2',e_2''$ are oriented outward at $w_2$.
If necessary we reflect the chart $\Gamma$,
we can assume that
the edge $e_4'$ is oriented from $w_4$ to $w_5$.
Since both of $e_2'',e_4'$ are oriented inward at $w_5$,
the edge $e_3''$ is oriented from $w_5$ to $w_3$.
Moreover, 
since both of $e_1'',e_3''$ are oriented inward at $w_3$
by (1),
the edge $e_3'$ is oriented from $w_3$ to $w_4$. 
Therefore $\Gamma$ contains the pseudo chart as shown
in Fig.~\ref{Fig23}(a).

{\bf Case (ii).}
One of $e_2',e_2''$ is oriented inward at $w_2$, and the other is oriented outward at $w_2$.
If necessary we reflect the chart $\Gamma$,
we can assume that 
the edge $e_2'$ is oriented inward at $w_2$, and the edge $e_2''$ is oriented outward at $w_2$.

Next, we shall show that $e_3'$ is oriented outward at $w_3$.
If $e_3'$ is oriented inward at $w_3$,
then $e_3''$ is oriented outward at $w_3$
(because, $e_1''$ is oriented inward at $w_3$ by (1)).
Thus both of $e_2''$ and $e_3''$ are oriented inward at $w_5$.
Hence the edge $e_4'$ is oriented from $w_5$ to $w_4$
(see Fig.~\ref{Fig23}(e)).
Thus, both of $\partial D_3$ and $\partial D_4$
are oriented clockwise or anticlockwise.
Hence,
by Lemma~\ref{Theorem3AngledDisk}, 
we have $w(\Gamma\cap{\rm Int}D_3)\geqq1$
and $w(\Gamma\cap{\rm Int}D_4)\geqq1$.
Moreover, since one of $D_1$ and $D_2$ is a 4-angled disk
with one feeler,
by Corollary~\ref{4AngledTwoFeeler}(a)
we have $w(\Gamma\cap{\rm Int}D_1)\geqq1$ or 
$w(\Gamma\cap{\rm Int}D_2)\geqq1$.
Thus\\
$\begin{array}{rl}
7=& w(\Gamma)\vspace{2mm}\\
= & w(G)+w(\Gamma\cap{\rm Int}D_1)+
w(\Gamma\cap{\rm Int}D_2)+w(\Gamma\cap{\rm Int}D_3)+
w(\Gamma\cap{\rm Int}D_4)\vspace{2mm}\\
\geqq & 5+1+1+1=8.\vspace{2mm}
\end{array}
$\\
This is a contradiction.
Hence $e_3'$ is oriented outward at $w_3$.

Since $e_1''$ is not middle at $w_3$,
the edge $e_3''$ is oriented inward at $w_3$.
If necessary we move the point $\infty$ in $D_3$,
we can assume that the edge $e_4'$ is oriented from 
$w_4$ to $w_5$.
Therefore, if $D_1$ (resp. $D_2$) has one feeler, then 
$\Gamma$ contains the pseudo chart as shown
in Fig.~\ref{Fig23}(b) (resp. Fig.~\ref{Fig23}(c)).
\end{Proof}


\begin{proposition}
\label{NoTypeF}
Let $\Gamma$ be a minimal chart of type $(m;5,2)$.
Then $\Gamma_m$ does not contain the graph as shown
in Fig.~\ref{Fig13}$($f$)$.
\end{proposition}

\begin{Proof}
Suppose that $\Gamma_m$ contains the graph as shown
in Fig.~\ref{Fig13}(f), say $G$.
Since the graph $G$ separates the 2-sphere $S^2$ 
into four disks.
Two of the four disks are 4-angled disks,
say $D_1,D_2$.
Two of the four disks are 3-angled disks,
say $D_3,D_4$.
Without loss of generality we can assume that
$D_4$ contains the point at infinity, $\infty$
(see Fig.~\ref{Fig23}(d)).
By Lemma~\ref{OriNoTypeF},
we can assume that
$\Gamma$ contains one of the three pseudo charts
as shown in Fig.~\ref{Fig23}(a),(b),(c).

Suppose that
$\Gamma$ contains  the pseudo chart
as shown in Fig.~\ref{Fig23}(a).
We use the notations as shown in Fig.~\ref{Fig23}(a),
where 
\begin{enumerate}
\item[(1)] 
$e_3,e_4,e_5$
are internal edges (possibly terminal edges)
of label $m+1$ oriented outward 
at $w_3,w_4,w_5$ in $D_4$,
respectively,
\end{enumerate}
but none of $e_3,e_4,e_5$
are middle at $w_3,w_4$ or $w_5$.
Thus by Assumption~\ref{AssumeTerminal},
\begin{enumerate}
\item[(2)] none of $e_3,e_4,e_5$ are terminal edges.
\end{enumerate}
Hence by IO-Calculation with respect to $\Gamma_{m+1}$
in $D_4$,
we have $w(\Gamma\cap{\rm Int}D_4)\geqq2$.
Moreover, since one of $D_1,D_2$ contains one feeler,
by Corollary~\ref{4AngledTwoFeeler}(a)
we have $w(\Gamma\cap{\rm Int}D_1)\geqq1$ or
$w(\Gamma\cap{\rm Int}D_2)\geqq1$.
Thus

$\begin{array}{rl}
7=& w(\Gamma)\vspace{2mm}\\
\geqq & w(G)+w(\Gamma\cap{\rm Int}D_1)+
w(\Gamma\cap{\rm Int}D_2)+
w(\Gamma\cap{\rm Int}D_4)\vspace{2mm}\\
\geqq & 5+1+2=8.\vspace{2mm}
\end{array}
$\\
This is a contradiction.
Hence $\Gamma$ does not contain the pseudo chart
as shown in Fig.~\ref{Fig23}(a).

Suppose that
$\Gamma$ contains the pseudo chart
as shown in Fig.~\ref{Fig23}(b).
Without loss of generality,
we can assume that $D_1$ has one feeler.
Thus by Corollary~\ref{4AngledTwoFeeler}(a)
we have $w(\Gamma\cap{\rm Int}D_1)\geqq1$.

We use the notations as shown in Fig.~\ref{Fig23}(b),
where
\begin{enumerate}
\item[(3)] 
${e}_1,{e}_2,{e}_5$
are internal edges (possibly terminal edges)
of label $m+1$ oriented outward 
at $w_1,w_2,w_5$ in $D_2$,
respectively,
\end{enumerate}
but neither ${e}_2$ nor ${e}_5$
is middle at $w_2$ or $w_5$.
Thus by Assumption~\ref{AssumeTerminal},
\begin{enumerate}
\item[(4)] neither ${e}_2$ nor ${e}_5$ 
is a terminal edge.
\end{enumerate}
Hence by IO-Caclulation  with respect to $\Gamma_{m+1}$
in $D_2$,
we have $w(\Gamma\cap{\rm Int}D_2)\geqq1$.

Since the boundary $\partial D_4$ is oriented anticlockwise,
by Lemma~\ref{Theorem3AngledDisk}
we have $w(\Gamma\cap{\rm Int}D_4)\geqq1$.
Thus 

$\begin{array}{rl}
7=& w(\Gamma)\vspace{2mm}\\
\geqq & w(G)+w(\Gamma\cap{\rm Int}D_1)+
w(\Gamma\cap{\rm Int}D_2)+
w(\Gamma\cap{\rm Int}D_4)\vspace{2mm}\\
\geqq & 5+1+1+1=8.\vspace{2mm}
\end{array}
$\\
This is a contradiction.
Hence $\Gamma$ does not contain the pseudo chart
as shown in Fig.~\ref{Fig23}(b).

Suppose that
$\Gamma$ contains  the pseudo chart
as shown in Fig.~\ref{Fig23}(c).
Without loss of generality,
we can assume that $D_2$ has one feeler.
Thus by Corollary~\ref{4AngledTwoFeeler}(a)
we have $w(\Gamma\cap{\rm Int}D_2)\geqq1$.

Since the boundary $\partial D_4$ is oriented anticlockwise,
by Lemma~\ref{Theorem3AngledDisk}
we have $w(\Gamma\cap{\rm Int}D_4)\geqq1$.
Hence the condition $w(\Gamma)=7$ implies that
\begin{enumerate}
\item[(5)]  $w(\Gamma\cap{\rm Int}D_1)=0$.
\end{enumerate}

We use the notations as shown in Fig.~\ref{Fig23}(c),
where ${e}_1,{e}_3,{e}_4$
are internal edges (possibly terminal edges)
of label $m+2$ at $w_1,w_3,w_4$ in $D_1$,
respectively, 
\begin{enumerate}
\item[(6)] 
${e}_1,{e}_3$ are
oriented outward at $w_1,w_3$, respectively,
\end{enumerate}
but neither ${e}_3$ nor ${e}_4$ 
is middle at $w_3$ or $w_4$.
Thus by Assumption~\ref{AssumeTerminal},
neither ${e}_3$  nor ${e}_4$ 
is a terminal edge.
Hence by (5) and (6),
we have $e_3=e_4$.
However there exists a lens.
This contradicts Lemma~\ref{NoLens}.
Thus $\Gamma$ does not contain the pseudo chart
as shown in Fig.~\ref{Fig23}(c).

Therefore we have a contradiction for all cases.
Hence $\Gamma_m$ does not contain the graph as shown
in Fig.~\ref{Fig13}(f).
\end{Proof}



\section{Triangle Lemma}
\label{s:TriangleLemma}

In this section,
we review Triangle Lemma.
These lemmas will be used in the next section.

\begin{figure}[htb]
\centerline{\includegraphics{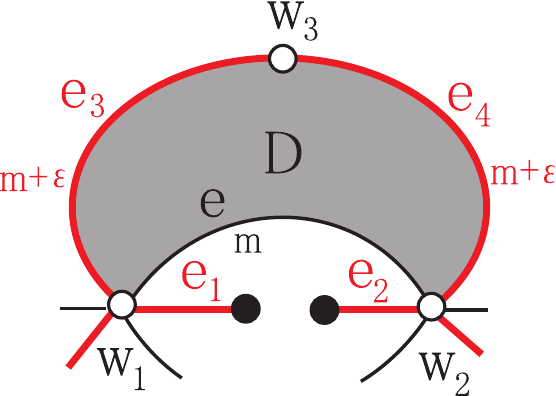}}
\caption{\label{Fig24}
The white vertices $w_1$ and $w_2$ are in $\Gamma_m\cap\Gamma_{m+\varepsilon}$, and $\varepsilon\in\{+1,-1\}$.}
\end{figure}

\begin{lemma}
\label{CorDiskLemma}
{\rm (\cite[Lemma 5.4]{ChartApp1})}
If a minimal chart $\Gamma$ contains the pseudo chart 
as shown in Fig.~\ref{Fig24}, 
then the interior of the disk $D$ contains at least one white vertex, 
where $D$ is the disk with the boundary $e_3\cup e_4\cup e$.
\end{lemma}

\begin{figure}[htb]
\centerline{\includegraphics{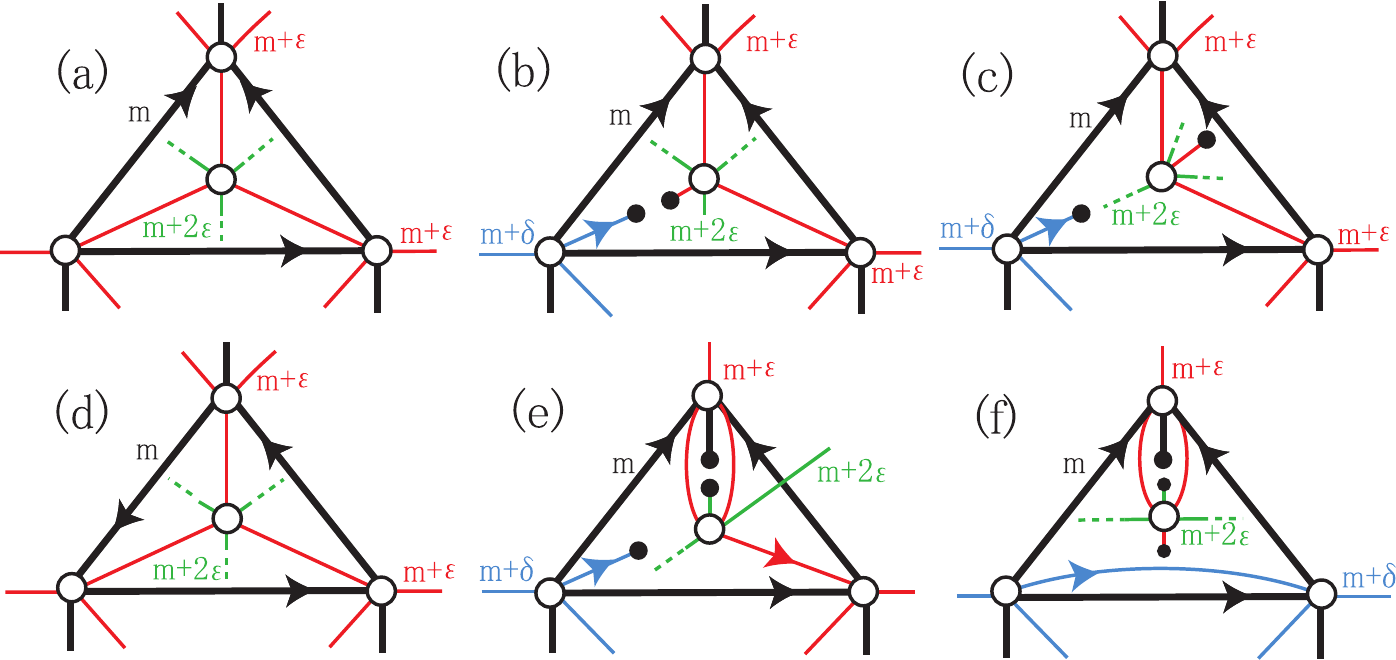}}
\caption{\label{Fig25}
(a),(b),(c),(d) 3-angled disks without feelers.
(e),(f) 3-angled disks with one feeler.}
\end{figure}

\begin{lemma}
{\rm (\cite[Lemma 4.2(b)]{ChartAppIX})}
\label{Theorem3AngledDiskChart9}
Let $\Gamma$ be a minimal chart, and $m$ a label of $\Gamma$.
Let $D$ be a special $3$-angled disk of $\Gamma_m$
with at most two feelers.
If $w(\Gamma\cap {\rm Int}D)=w(\Gamma_{m+\varepsilon}\cap {\rm Int}D)=1$ for some $\varepsilon\in\{+1,-1\}$,
then a regular neighborhood of $D$ contains one of the RO-families of the six pseudo charts as shown in 
Fig.~\ref{Fig25}.
\end{lemma}

\begin{figure}
\centerline{\includegraphics{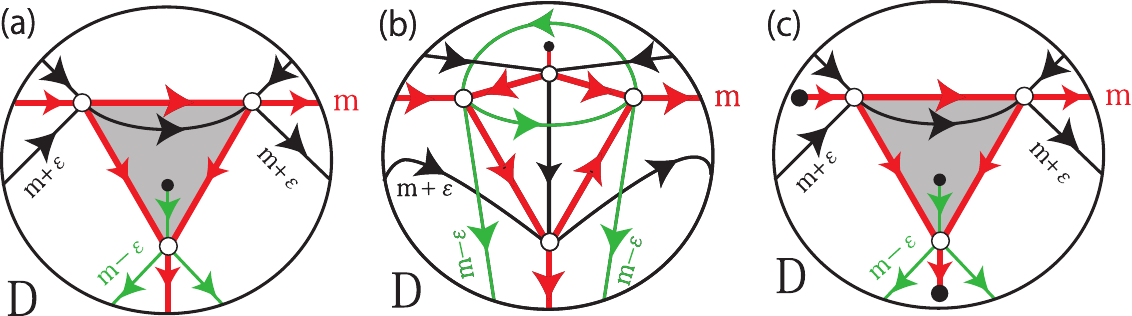}}
\caption{\label{Fig26}
The gray region is the 3-angled disk $D_1$. 
The thick lines are edges of label $m$,
and $\varepsilon\in\{+1,-1\}$.}
\end{figure}

\begin{lemma}$($Triangle Lemma$)$
{\rm (\cite[Lemma 8.3]{ChartAppIV})}
\label{LemmaTriangle}
\begin{enumerate}
\item[{\rm (a)}]
For a chart $\Gamma$,
if there exists  a $3$-angled disk $D_1$ of $\Gamma_m$ without feelers in a disk $D$ as shown in Fig.~\ref{Fig26}$($a$)$ and if $w(\Gamma\cap${\rm Int}$D_1)=0$,
then there exists a chart obtained from $\Gamma$ by C-moves in $D$ which contains the pseudo chart in $D$ as shown in Fig.~\ref{Fig26}$($b$)$. 
\item[{\rm (b)}] For a minimal chart $\Gamma$, 
if there exists a $3$-angled disk $D_1$ of $\Gamma_m$ without feelers in a disk $D$ as shown in Fig.~\ref{Fig26}$($c$)$,
then $w(\Gamma\cap${\rm Int}$D_1)\ge1$.
\end{enumerate}
\end{lemma}

\begin{figure}
\centerline{\includegraphics{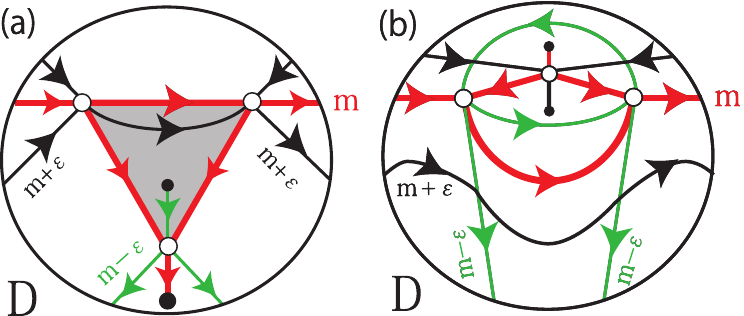}}
\caption{\label{Fig27}
The gray region is the 3-angled disk $D_1$. 
The thick lines are edges of label $m$,
and $\varepsilon\in\{+1,-1\}$.}
\end{figure}

By the above lemma,
we can show the following corollary
by using C-II moves and a C-III move:

\begin{corollary}
\label{CorTriangle}
For a chart $\Gamma$,
if there exists  a $3$-angled disk $D_1$ of $\Gamma_m$ without feelers in a disk $D$ as shown in Fig.~\ref{Fig27}$($a$)$ and if $w(\Gamma\cap${\rm Int}$D_1)=0$,
then there exists a chart obtained from $\Gamma$ by C-moves in $D$ which contains the pseudo chart in $D$ as shown in Fig.~\ref{Fig27}$($b$)$. 
\end{corollary}

\begin{lemma}
{\rm (\cite[Theorem 1.1]{ChartAppIX})}
\label{NoMinimal(43)} 
There is 
no minimal chart of 
type $(4,3)$.
\end{lemma}



\section{Case of the graph as shown in Fig.~\ref{Fig13}(g)}
\label{s:TypeI}

In this section,
we shall show that if $\Gamma$ is a minimal chart of type $(m;5,2)$,
then the graph $\Gamma_m$ does not contain
the graph as shown in Fig.~\ref{Fig13}(g).
Moreover, we shall show the main theorem.

\begin{figure}[htb]
\centerline{\includegraphics{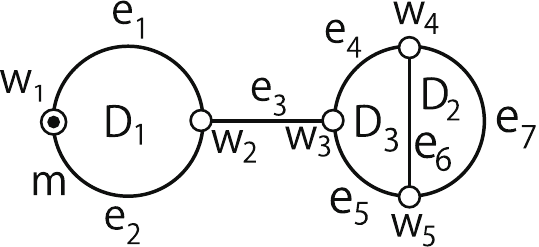}}
\caption{\label{Fig28}
The graph as shown in Fig.~\ref{Fig13}(g).}
\end{figure}

Suppose that 
the graph $\Gamma_m$ contains the graph as shown
in Fig.~\ref{Fig13}$($g$)$.
Form now on throughout this section,
we use the notations as shown in Fig.~\ref{Fig28},
where
\begin{enumerate}
\item[(a)] $w_1,w_2,\cdots,w_5$ are five white vertices,
and 
\item[(b)] 
$e_1,e_2,\cdots,e_7$ are seven internal edges of label $m$
with
$\partial e_1=\partial e_2=\{w_1,w_2\}$,
$\partial e_3=\{w_2,w_3\}$,
$\partial e_4=\{w_3,w_4\}$,
$\partial e_5=\{w_3,w_5\}$,
$\partial e_6=\partial e_7=\{w_4,w_5\}$,
\item[(c)] $D_1,D_2$ are special 2-angled disks
with $\partial D_1=e_1\cup e_2$ and 
$\partial D_2=e_6\cup e_7$,
\item[(d)] $D_3$ is the special 3-angled disk
with $\partial D_3=e_4\cup e_5\cup e_6$.
\end{enumerate}

\begin{figure}[htb]
\centerline{\includegraphics{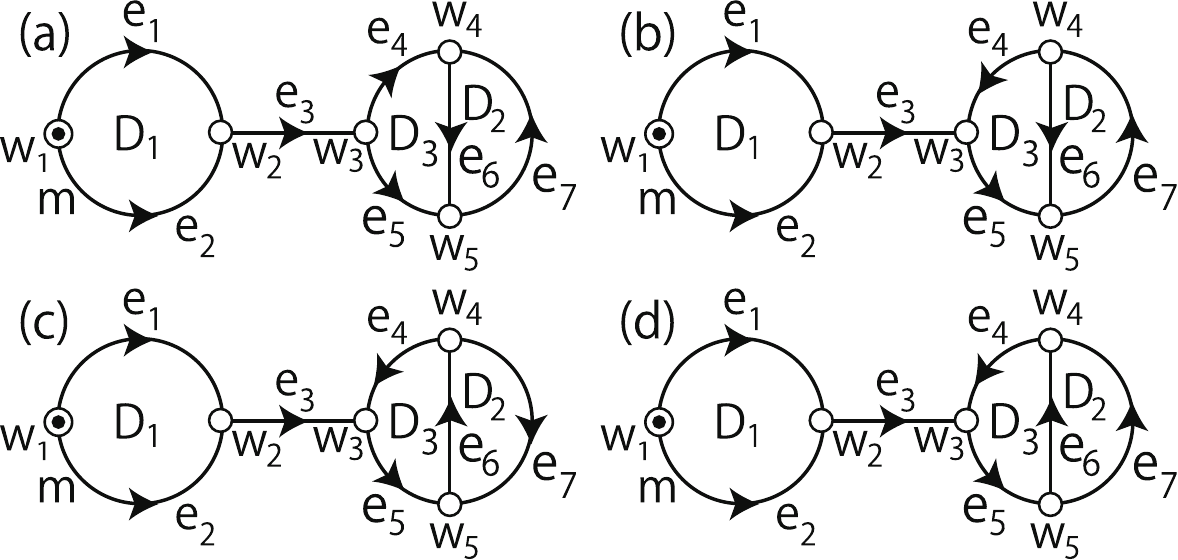}}
\caption{\label{Fig29}
The graphs as shown in Fig.~\ref{Fig13}(g).}
\end{figure}

\begin{lemma}
\label{OriNoTypeI}
Let $\Gamma$ be a minimal chart of type $(m;5,2)$.
If $\Gamma_m$ contains the graph as shown
in Fig.~\ref{Fig13}$($g$)$,
then $\Gamma_m$ contains one of RO-families
of the four graphs as shown in 
Fig.~\ref{Fig29}.
\end{lemma}

\begin{Proof}
We use the notations as shown in Fig.~\ref{Fig28}.
Without loss of generality,
we can assume that the terminal edge of label $m$ at $w_1$
is oriented inward at $w_1$.
Then by Assumption~\ref{AssumeTerminal},
both of $e_1,e_2$ are oriented from $w_1$ to $w_2$.
Thus the edge $e_3$ is oriented form $w_2$ to $w_3$.
There are two cases:
(i) $e_3$ is middle at $w_3$,
(ii) $e_3$ is not middle at $w_3$.

{\bf Case (i).}
Since $e_3$ is middle at $w_3$,
both of $e_4,e_5$ are oriented outward at $w_3$.
Hence 
\begin{enumerate}
\item[(1)] $e_5$ is oriented inward at $w_5$. 
\end{enumerate}

If necessary we reflect the chart $\Gamma$,
we can assume that $e_6$ is oriented from $w_4$ to $w_5$.
Thus by (1),
the edge $e_7$ is oriented from $w_5$ to $w_4$.
Hence  $\Gamma_m$ contains the graph as shown in
Fig.~\ref{Fig29}(a).

{\bf Case (ii).}
Since $e_3$ is not middle at $w_3$,
one of $e_4,e_5$ is oriented inward at $w_3$
and the other is oriented outward at $w_3$.
If necessary we reflect the chart $\Gamma$,
we can assume that
$e_4$ is oriented inward at $w_3$
and $e_5$ is oriented outward at $w_3$.
Thus 
\begin{enumerate}
\item[(2)] $e_5$ is oriented inward at $w_5$. 
\end{enumerate}

If $e_6$ is oriented from $w_4$ to $w_5$,
then by (2) the edge $e_7$ is oriented from $w_5$ to $w_4$.
Hence $\Gamma_m$ contains the graph as shown in
Fig.~\ref{Fig29}(b).

If $e_6$ is oriented from $w_5$ to $w_4$,
then $\Gamma_m$ contains one of the two graphs as shown in
Fig.~\ref{Fig29}(c),(d).
\end{Proof}


\begin{lemma}
\label{NoTypeICaesC}
Let $\Gamma$ be a minimal chart of type $(m;5,2)$.
Then $\Gamma_m$ does not contain the graph as shown
in Fig.~\ref{Fig29}$($c$)$.
\end{lemma}

\begin{Proof}
Suppose that
$\Gamma_m$ contains the graph as shown
in Fig.~\ref{Fig29}(c), say $G$.
We use the notations as shown in Fig.~\ref{Fig28}
and Fig.~\ref{Fig29}(c).

Since $\partial D_2$ is oriented clockwise,
since $e_4$ is oriented outward at $w_4$ and
since $e_5$ is oriented inward at $w_5$,
by Lemma~\ref{Theorem2AngledDisk2} we have
 $w(\Gamma\cap{\rm Int}D_2)\geqq2$. 

Since $\partial D_3$ is oriented anticlockwise,
by Lemma~\ref{Theorem3AngledDisk}
we have $w(\Gamma\cap{\rm Int}D_3)\geqq1$. 
Thus we have\vspace{2mm}

$\begin{array}{rcl}
7=w(\Gamma) & \geqq & w(G)+w(\Gamma\cap{\rm Int}D_2)+
w(\Gamma\cap{\rm Int}D_3)\vspace{2mm}\\
& \geqq & 5+2+1=8.\vspace{2mm}
\end{array}
$\\
This is a contradiction.
Therefore $\Gamma_m$ does not contain the graph as shown
in Fig.~\ref{Fig29}(c).
\end{Proof}

\begin{lemma}
\label{NoTypeICaesB}
Let $\Gamma$ be a minimal chart of type $(m;5,2)$.
Then $\Gamma_m$ does not contain the graph as shown
in Fig.~\ref{Fig29}$($b$)$.
\end{lemma}

\begin{figure}[htb]
\centerline{\includegraphics{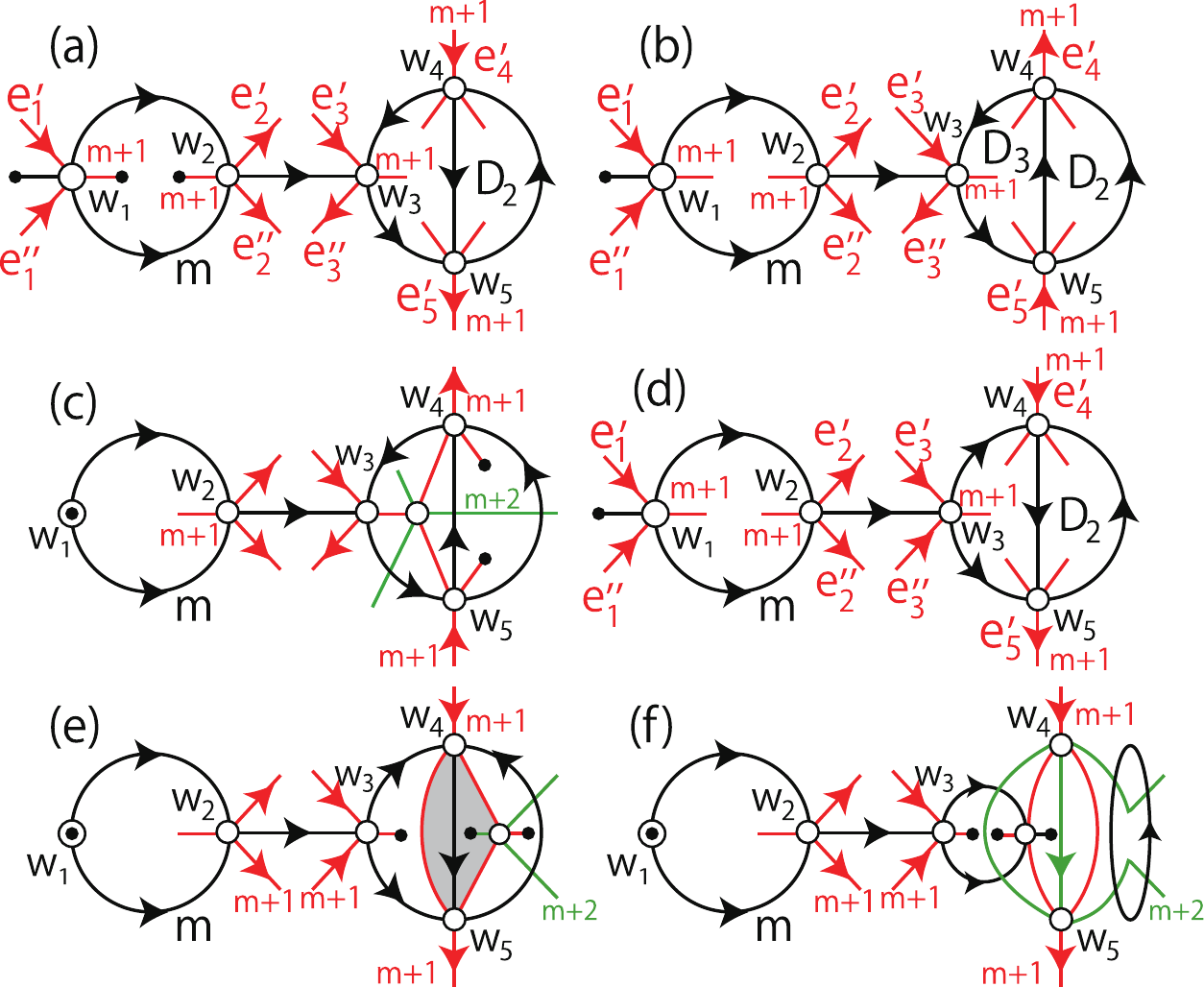}}
\caption{\label{Fig30}
The graphs as shown in Fig.~\ref{Fig13}(g).
The gray region is the 3-angled disk $D$.}
\end{figure}

\begin{Proof}
Suppose that
$\Gamma_m$ contains the graph as shown
in Fig.~\ref{Fig29}(b).
We use the notations as shown in Fig.~\ref{Fig28}
and Fig.~\ref{Fig29}(b).

By the similar way of the proof of Lemma~\ref{NoTypeICaesC},
we have  $w(\Gamma\cap{\rm Int}D_2)\geqq2$. 
Thus the condition $w(\Gamma)=7$ implies that
\begin{enumerate}
\item[(1)] $w(\Gamma\cap{\rm Int}D_1)=0$ and
$w(\Gamma\cap(S^2-(D_1\cup D_2\cup D_3)))=0$. 
\end{enumerate}
Hence by Lemma~\ref{Theorem2AngledDisk},
a regular neighborhood of $D_1$
contains the pseudo chart as shown in Fig.~\ref{Fig10}(b)
(see Fig.~\ref{Fig30}(a)).

We use the notations as shown in Fig.~\ref{Fig30}(a),
where 
\begin{enumerate}
\item[(2)] $e_1',e_1'',e_3',e_4'$
are internal edges (possibly terminal edges)
of label $m+1$ oriented inward at $w_1,w_1,w_3,w_4$,
respectively,
\item[(3)] $e_2',e_2'',e_3'',e_5'$
are internal edges (possibly terminal edges)
of label $m+1$ oriented outward at $w_2,w_2,w_3,w_5$,
respectively.
\end{enumerate}
Moreover, none of  $e_2',e_2'',e_3'',e_5'$ are
middle at $w_2,w_3$ or $w_5$.
Thus by Assumption~\ref{AssumeTerminal},
\begin{enumerate}
\item[(4)] none of $e_2',e_2'',e_3'',e_5'$ are
terminal edges.
\end{enumerate}
Hence by (1),(2),(3),
none of $e_1',e_1'',e_3',e_4'$ are
terminal edges.
Thus for the edge $e_2''$,
we have $e_2''=e_1''$.
However $e_2\cup e_2''$ bounds a lens.
This contradicts Lemma~\ref{NoLens}.
Therefore $\Gamma_m$ does not contain the graph as shown
in Fig.~\ref{Fig29}(b).
\end{Proof}

\begin{lemma}
\label{NoTypeICaesD}
Let $\Gamma$ be a minimal chart of type $(m;5,2)$.
Then $\Gamma_m$ does not contain the graph as shown
in Fig.~\ref{Fig29}$($d$)$.
\end{lemma}

\begin{Proof}
Suppose that
$\Gamma_m$ contains the graph as shown
in Fig.~\ref{Fig29}(d).
We use the notations as shown in Fig.~\ref{Fig28}
and Fig.~\ref{Fig29}(d).

Since $\partial D_3$ is oriented anticlockwise,
by Lemma~\ref{Theorem3AngledDisk}
we have 
\begin{enumerate}
\item[(1)] $w(\Gamma\cap{\rm Int}D_3)\geqq1$. 
\end{enumerate}

{\bf Claim.} $w(\Gamma\cap{\rm Int}D_2)=0$ and
$w(\Gamma\cap{\rm Int}D_3)=1$. 

{\it Proof of Claim~$1$.}
Let $e$ be the terminal edge of label $m$ at $w_1$.

If $e\subset D_1$, then by Lemma~\ref{Theorem2AngledDisk}
we have $w(\Gamma\cap{\rm Int}D_1)\geqq1$.
Thus by (1) and $w(\Gamma)=7$,
we have $w(\Gamma\cap{\rm Int}D_2)=0$
and $w(\Gamma\cap{\rm Int}D_3)=1$.

Now, suppose  $e\not\subset D_1$.
We use the notations as shown in Fig.~\ref{Fig30}(b),
where 
\begin{enumerate}
\item[(2)] $e_1',e_1'',e_3',e_5'$ are
internal edges (possibly terminal edges) of label $m+1$
oriented inward at $w_1,w_1,w_3,w_5$, 
respectively, and
\item[(3)] $e_2',e_2'',e_3'',e_4'$ are
internal edges (possibly terminal edges) of label $m+1$
oriented outward at $w_2,w_2,w_3,w_4$, 
respectively.
\end{enumerate}
Moreover, none of $e_2',e_2'',e_3'',e_4'$
are middle at $w_2,w_3$ or $w_4$.
Thus by Assumption~\ref{AssumeTerminal},
none of the four edges $e_2',e_2'',e_3'',e_4'$
are terminal edges.

If $w(\Gamma\cap(S^2-(D_1\cup D_2\cup D_3)))=0$,
then by (2) and (3)
none of the four edges $e_1',e_1'',e_3',e_5'$ are
terminal edges.
Hence for the edge $e_2''$,
we have $e_2''=e_1''$.
However, there exists a lens.
This contradicts Lemma~\ref{NoLens}.
Thus $w(\Gamma\cap(S^2-(D_1\cup D_2\cup D_3)))\geqq1$.

Hence by (1) and $w(\Gamma)=7$,
we have $w(\Gamma\cap{\rm Int}D_2)=0$
and $w(\Gamma\cap{\rm Int}D_3)=1$.
Thus Claim holds. \hfill {$\square$}\vspace{1.5em}

Since $\partial D_3$ is oriented anticlockwise,
by Claim and Lemma~\ref{Theorem3AngledDiskChart9}
a regular neighborhood of $D_3$ contains
the pseudo chart as shown in Fig.~\ref{Fig25}(d).
Moreover,
by Claim and Lemma~\ref{Theorem2AngledDisk}
a regular neighborhood of $D_2$ contains
the pseudo chart as shown in Fig.~\ref{Fig10}(b)
(see Fig.~\ref{Fig30}(c)).
Hence the chart $\Gamma$
contains the pseudo chart as shown in Fig.~\ref{Fig24}.
Thus by Lemma~\ref{CorDiskLemma},
we have $w(\Gamma\cap {\rm Int}D_3)\ge2$.
This contradicts Claim.
Thus we complete the proof of Lemma~\ref{NoTypeICaesD}.
\end{Proof}

\begin{proposition}
\label{NoTypeI}
Let $\Gamma$ be a minimal chart of type $(m;5,2)$.
Then $\Gamma_m$ does not contain the graph as shown
in Fig.~\ref{Fig13}$($g$)$.
\end{proposition}

\begin{Proof}
Suppose that
$\Gamma_m$ contains the graph as shown in 
Fig.~\ref{Fig13}(g).
Then by Lemma~\ref{OriNoTypeI},
the graph $\Gamma_m$ contains
one of RO-families of the four graphs as shown in
Fig.~\ref{Fig29}.
Hence by Lemma~\ref{NoTypeICaesC},
Lemma~\ref{NoTypeICaesB} and
Lemma~\ref{NoTypeICaesD},
the graph $\Gamma_m$ contains
one of the RO-family of the graph as shown in
Fig.~\ref{Fig29}(a).
Without loss of generality,
we can assume that the graph $\Gamma_m$ contains
 of the graph as shown in
Fig.~\ref{Fig29}(a).

Since $\partial D_2$ is oriented anticlockwise,
by Lemma~\ref{Theorem2AngledDisk}
we have 
\begin{enumerate}
\item[(1)] $w(\Gamma\cap{\rm Int}D_2)\geqq1$. 
\end{enumerate}

{\bf Claim.} $w(\Gamma\cap{\rm Int}D_2)=1$ and $w(\Gamma\cap{\rm Int}D_3)=0$.

{\it Proof of Claim.}
Let $e$ be the terminal edge of label $m$ at $w_1$.

If $e\subset D_1$, then by Lemma~\ref{Theorem2AngledDisk}
we have $w(\Gamma\cap{\rm Int}D_1)\geqq1$.
Thus by (1) and $w(\Gamma)=7$,
we have $w(\Gamma\cap{\rm Int}D_2)=1$
and $w(\Gamma\cap{\rm Int}D_3)=0$.

Now, suppose $e\not\subset D_1$.
We use the notations as shown in Fig.~\ref{Fig30}(d), 
where 
\begin{enumerate}
\item[(2)] $e_1',e_1'',e_3',e_3'',e_4'$
are internal edges (possibly terminal edges)
of label $m+1$ oriented inward at $w_1,w_1,w_3,w_3,w_4$,
respectively.
\end{enumerate}
Moreover, none of $e_1',e_1'',e_3',e_3''$ are
middle at $w_1$ or $w_3$.
Thus by Assumption~\ref{AssumeTerminal},
none of the four edges
$e_1',e_1'',e_3',e_3''$ are terminal edges.
Thus by (2) and by
IO-Calculation with respect to $\Gamma_{m+1}$ in 
$Cl(S^2-(D_1\cup D_2\cup D_3))$,
we have $w(\Gamma\cap(S^2-(D_1\cup D_2\cup D_3)))\geqq1$.
Hence by (1) and $w(\Gamma)=7$,
we have $w(\Gamma\cap{\rm Int}D_2)=1$ and
$w(\Gamma\cap{\rm Int}D_3)=0$.
Thus Claim holds. \hfill {$\square$}\vspace{1.5em}

By Claim and Lemma~\ref{Theorem3AngledDisk},
a regular neighborhood of $D_3$
contains the pseudo chart as shown in Fig.~\ref{Fig16}(b).
Moreover, by Claim and Lemma~\ref{Theorem2AngledDisk2},
a regular neighborhood of $D_2$
contains one of the two pseudo charts as shown 
in Fig.~\ref{Fig21}(b),(c).
Hence there exists a 3-angled disk $D$ of $\Gamma_{m+1}$
in $D_2\cup D_3$.

Let $w_6$ be the white vertex in ${\rm Int}D_2$,
and $e'$ the terminal edge of label $m+1$ at $w_6$.
Since $w(\Gamma\cap {\rm Int}D)=0$ 
by Claim,
a regular neighborhood of $D$ contains the pseudo chart as shown in Fig.~\ref{Fig16}(a).
Hence $e'\not\subset D$
(see Fig.~\ref{Fig30}(e)).
Thus $\Gamma$ contains the pseudo chart as shown in
Fig.~\ref{Fig27}(a).
Thus by Corollary~\ref{CorTriangle},
there exists a minimal chart $\Gamma'$ obtained from $\Gamma$
by C-moves which contains the pseudo chart
as shown in Fig.~\ref{Fig27}(b)
(see Fig.~\ref{Fig30}(f)).
Hence $\Gamma$ is C-move equivalnet to 
the minimal chart $\Gamma'$ of type $(m;4,3)$.
This contradicts Lemma~\ref{NoMinimal(43)}.
Hence $\Gamma_m$ does not contain the graph as shown
in Fig.~\ref{Fig29}(a).

Hence $\Gamma_m$ does not contains the graph
as shown in Fig.~\ref{Fig13}(g).
Therefore we complete the proof of Proposition~\ref{NoTypeI}.
\end{Proof}

\begin{lemma}$(${\rm \cite[Theorem 1.1]{ChartAppIV}}$)$
\label{LemmaNoLoop}
There is no loop in any minimal chart with exactly seven white vertices.
\end{lemma}

Now, we shall show the main theorem.

{\it Proof of Theorem~\ref{MainTheorem}.}
Let $\Gamma$ be a minimal chart of type $(m;5,2)$.
Suppose that there exists a connected component $G$ 
of $\Gamma_m$ with $w(G)=5$.
Then by Lemma~\ref{LemmaNoLoop},
the graph $G$ does not contain any loop.
Thus by Lemma~\ref{LemmaWithTerminal},
the graph $G$ is one of nine graphs as shown 
in  Fig.~\ref{Fig02} and Fig.~\ref{Fig13}.
Hence the main theorem follows from the seven propositions
(Lemma~\ref{NoTypeA}, Lemma~\ref{NoTypeC},
Lemma~\ref{NoTypeB}, Lemma~\ref{NoTypeD},
Lemma~\ref{NoTypeE}, Proposition~\ref{NoTypeF}
and Proposition~\ref{NoTypeI}).
Therefore we complete the proof of the main theorem.
\hfill {$\square$}\vspace{1.5em}




\vspace{5mm}

\begin{minipage}{65mm}
{Teruo NAGASE
\\
{\small Tokai University \\
4-1-1 Kitakaname, Hiratuka \\
Kanagawa, 259-1292 Japan\\
\\
nagase@keyaki.cc.u-tokai.ac.jp
}}
\end{minipage}
\begin{minipage}{65mm}
{Akiko SHIMA 
\\
{\small Department of Mathematics, 
\\
Tokai University
\\
4-1-1 Kitakaname, Hiratuka \\
Kanagawa, 259-1292 Japan\\
shima@keyaki.cc.u-tokai.ac.jp\\
shima-a@tokai.ac.jp
}}
\end{minipage}


\vspace{0.7cm}

{\bf List of terminologies}\vspace{2mm}\\
{\small $
\begin{array}{ll||}
\text{$k$-angled disk} & p9 \\
\text{BW-vertex} & p2 \\
\text{C-move equivalent} & p4 \\
\text{chart} & p3 \\
\text{complexity $(w(\Gamma),-f(\Gamma))$} & p4 \\
\text{feeler} & p10 \\
\text{free edge} & p4 \\
\text{hoop} & p4 \\
\text{internal edge} & p6 \\
\text{inward} & p4 \\
\text{inward arc} & p22 \\
\text{IO-Calculation} & p23 \\
\text{keeping $X$ fixed} & p24 \\
\text{lens} & p6 \\
\text{locally minimal} & p25 \\
\end{array}
~~
\begin{array}{ll}
\text{loop} & p6 \\
\text{middle arc} & p4 \\
\text{middle at $v$} & p4 \\
\text{minimal chart} & p4 \\
\text{outward} & p4 \\
\text{outward arc} & p22 \\
\text{point at infinity $\infty$} & p4 \\
\text{pseudo chart} & p12 \\
\text{ring} & p4 \\
\text{RO-family} & p13 \\
\text{simple hoop} & p5 \\
\text{special $k$-angled disk} & p10 \\
\text{terminal edge} & p2 \\
\text{type $(m;n_1,n_2,\cdots,n_k)$ for a chart} & p2 \\
& \\
\end{array}
$}

\vspace{0.5cm}

{\bf List of notations}\vspace{2mm}\\
{\small $
\begin{array}{ll||}
\text{$\Gamma_m$} & p2 \\
\text{$w(\Gamma)$} & p4 \\
\text{$f(\Gamma)$} & p4 \\
\text{${\rm Int}X$} & p5 \\
\text{$\partial X$} & p5 \\
\text{$Cl(X)$} & p5 \\
\end{array}
$~~
$\begin{array}{ll}
\text{$\partial \alpha$} & p5 \\
\text{${\rm Int}\alpha$} & p5\\
\text{$a_{ij},b_{ij}$} & p6 \\
\text{$w(X)$} & p8 \\
\text{$c(X)$} & p24\\
& \\
\end{array}
$
}

\end{document}